\documentclass[11pt, reqno]{amsart}
\usepackage{amsmath,amsfonts,amssymb,amsmath,latexsym,mathrsfs}
\usepackage[all]{xy}
\usepackage{setspace}

\newtheorem{theorem}{Th\'eor\`eme}[section]
\newtheorem{lemma}[theorem]{Lemme}
\newtheorem{proposition}[theorem]{Proposition}
 
\newtheorem{definition}[theorem]{D\'efinition}
\theoremstyle{definition}  
\newtheorem{example}[theorem]{Example}

\newtheorem{remark}[theorem]{Remarque}

\newcommand{\HH}{{\mathfrak{H}}}\newcommand{\on}{\operatorname}

\renewcommand{\i}{\on{i}}\newcommand{\ul}{\underline}

\newcommand{\CC}{{\mathbb{C}}}\newcommand{\QQ}{{\mathbb{Q}}}

\newcommand{\ZZ}{{\mathbb{Z}}}

\newcommand{\f}{{\mathfrak{f}}}

\renewcommand{\t}{{\mathfrak{t}}}

\makeatother
\input cyracc.def
\DeclareFontFamily{U}{russian}{}
\DeclareFontShape{U}{russian}{m}{n}
        { <5><6> wncyr5
        <7><8><9> wncyr7
        <10><10.95><12><14.4><17.28><20.74><24.88> wncyr10 }{}
\DeclareSymbolFont{Russian}{U}{russian}{m}{n}
\DeclareSymbolFontAlphabet{\mathcyr}{Russian}
\makeatletter
\let\@math@cyr\mathcyr
\renewcommand{\mathcyr}[1]{\@math@cyr{\cyracc #1}}
\makeatother
\newcommand{\sh}{{\mathcyr{sh}}} 

\begin{document}

\title[Analogues elliptiques des nombres multiz\'etas]{Analogues elliptiques des 
nombres multiz\'etas}

\renewcommand{\abstractname}{R\' esum\' e}
\begin{abstract} 
Nous \'etudions des fonctions du param\`etre elliptique d\'efinies commes 
int\'egrales it\'er\'ees de fonctions elliptiques. Nous \'etablissons leur lien 
avec les ``associateurs elliptiques'' de notre pr\'ec\'edent travail au moyen 
de r\'ealisations fonctionnelles d'alg\`ebres de Lie apparaissant dans cette th\'eorie. 

\end{abstract}   

\author{Benjamin Enriquez}
\address{IRMA (CNRS), Universit\'e de Strasbourg, 7 
rue Ren\'e Descartes, F-67084 Strasbourg, France}
\email{b.enriquez@math.unistra.fr}

\maketitle

\section*{Introduction}

La th\'eorie des nombres multiz\'etas a d\'ebut\'e dans \cite{Za} avec la
construction de familles de relations entre ces nombres, reposant en partie sur 
le lien observ\'e par 
Kontsevich entre multiz\'etas et int\'egrales it\'er\'ees sur les espaces de modules 
de courbes rationnelles avec points marqu\'es. Parall\`element, 
Drinfeld a \'etabli des relations d'origine g\'eom\'etrique satisfaites par une 
s\'erie non-commutative, l'associateur KZ (\cite{Dr}) ; Le et Murakami ont 
identifi\'e l'associateur KZ \`a une s\'erie g\'en\'eratrice des multiz\'etas
(\cite{LM}), ce qui a permis de consid\'erer les relations de l'associateur KZ
comme un deuxi\`eme syst\`eme de relations entre multiz\'etas. 
Le lien entre les deux syst\`emes de relations a \'et\'e \'etudi\'e
par Furusho (\cite{Fu}). 

Un analogue elliptique de la th\'eorie des associateurs a \'et\'e construit dans 
\cite{En2} \`a partir d'un analogue elliptique de la connection de 
Knizhnik-Zamolodchikov (\cite{CEE}, voir aussi \cite{LR}). Le r\^{o}le de 
l'associateur KZ y est tenu par un couple de fonctions
$(A(\tau),B(\tau))$ d'un param\`etre $\tau$ dans le demi-plan de Poincar\'e, \`a valeurs 
dans un groupe de s\'eries non-commutatives \`a deux variables 
$\on{exp}(\hat\f_{2})$. Les r\'esultats principaux de \cite{En2}
sur le couple $(A(\tau),B(\tau))$ sont : le comportement de ce couple sous les transformations 
modulaires ; une famille de relations alg\'ebriques (d'origine g\'eom\'etrique)
satisfaites par $(A(\tau),B(\tau))$ ;  une \'equation diff\'erentielle satisfaite 
par le m\^{e}me objet ; son comportement en $\tau\to\on{i}\infty$. 
Un corollaire de cette \'etude est une famille de relations
alg\'ebriques entre int\'egrales it\'er\'ees de s\'eries d'Eisenstein et multiz\'etas. 
Un r\^{o}le important est jou\'e dans cette th\'eorie, et \'egalement dans la th\'eorie 
reli\'ee des motifs elliptiques universels (\cite{HM,Pk}), par une alg\`ebre de Lie 
$\langle \delta_{2n},n\geq -1\rangle \subset \on{Der}(\f_{2})$.  Nous rappelons ces r\'esultats en section 1. 

Le but principal de cet article est l'\'etude des coefficients des s\'eries $A(\tau),B(\tau)$.
Il s'agit de fonctions 
$$
I_{\ul{d}}(\tau), J_{\ul{d}}(\tau), \quad \ul d = (d_{1},\ldots,d_{n})\in 
\{-1,0,1,\ldots\}^{n}
$$
du param\`etre elliptique, qui sont des analogues elliptiques des nombres multiz\'etas. 

La section 2 est consacr\'ee \`a la d\'etermination d'expressions int\'egrales pour ces fonctions. Nous utilisons pour cela le calcul 
de l'holonomie r\'egularis\'ee des \'equations diff\'erentielles sur $]0,1[$ \`a valeurs dans une alg\`ebre libre, avec singularit\'es aux 
extr\'emit\'es. Ce calcul a \'et\'e effectu\'e 
dans \cite{En1} \`a partir d'id\'ees contenues dans \cite{LM}. Le r\'esultat de \cite{En1} est formul\'e en sections 2.1 et 2.2, et appliqu\'e 
en sections 2.3 et 2.4 au calcul d'expressions int\'egrales pour les $I_{\underline d}(\tau)$,  $J_{\underline d}(\tau)$ (relations (\ref{def:I}), 
(\ref{def:J}), (\ref{explicit:I}), (\ref{condensed:I})). 
En section 2.5, nous traduisons en termes des $I_{\underline d}(\tau)$,  $J_{\underline d}(\tau)$ certaines identit\'es satisfaites 
par $(A(\tau),B(\tau))$. 

La section 3 est consacr\'ee aux syst\`emes diff\'erentiels satisfaits par les $I_{\underline d}(\tau)$, $J_{\underline d}(\tau)$. En section 3.1, 
on construit des syst\`emes diff\'erentiels satisfaits par des int\'egrales it\'er\'ees g\'en\'erales du type de celles introduites en section 2.2. 
En section 3.2, on applique ce r\'esultat aux fonctions $I_{\underline d}(\tau)$,  $J_{\underline d}(\tau)$ et on obtient ainsi un syst\`eme 
diff\'erentiel satisfait par ces fonctions (th\'eor\`eme \ref{thm:ode}).  

En section 4, on montre l'\'equivalence entre ce syst\`eme diff\'erentiel et le syst\`eme diff\'erentiel satisfait par $(A(\tau),B(\tau))$
explicit\'e dans \cite{En2}. Ceci est r\'ealis\'e au moyen d'une r\'ealisation fonctionnelle de l'alg\`ebre de Lie 
$\langle\delta_{2n},n\geq -1\rangle$ (section \ref{sect:real:fonct}).

Enfin, en section \ref{section:DA}, on applique le syst\`eme diff\'erentiel du th\'eor\`eme \ref{thm:ode} au d\'eveloppement 
asymptotique des $I_{\ul{d}}(\tau), J_{\ul{d}}(\tau)$ en $\tau\to\on{i}\infty$ ; on montre que ce d\'eveloppement 
asymptotique s'exprime \`a l'aide de nombres multiz\'etas.

Signalons enfin les liens possibles entre le pr\'esent travail et \cite{BL} : les auteurs de cet article construisent une th\'eorie 
des polylogarithmes elliptiques multiples, qui sont certaines fonctions multivalu\'ees sur la vari\'et\'e $E_\tau^n$, o\`u 
$E_\tau :={\mathbb{C}}/({\mathbb{Z}} + \tau {\mathbb{Z}})$.  Ils projettent d'en d\'eduire, par sp\'ecialisation, des fonctions 
de $\tau$ qu'ils appellent ``fonctions multiz\'etas elliptiques''. On peut s'attendre \`a ce que ces fonctions pr\'esentent des liens 
\'etroits avec les fonctions $I_{\underline{d}}(\tau)$, $J_{\underline{d}}(\tau)$ du pr\'esent article. 

\section{Pr\'eliminaires : associateurs elliptiques}

Dans cette section, nous rappelons la construction et les propri\'et\'es de la fonction 
$\tau\mapsto (A(\tau),B(\tau))$ (\cite{En2}). 

\subsection{D\'efinition de $(A(\tau),B(\tau))$}\label{sec:def:A:B}

Soit pour $n\geq 2$,  $\bar\t_{1,n}$ l'alg\`ebre de Lie pr\'esent\'ee par les 
g\'en\'erateurs $x_{i},y_{i}$, $i\in\{1,\ldots,n\}$ et les relations 
$\sum_{i}x_{i} = \sum_{i}y_{i}=0$, $[x_{i},x_{j}] = [y_{i},y_{j}] = 0$, 
$[x_{i},y_{j}] = [x_{j},y_{i}]=:t_{ij}$ si $i\neq j$, $[x_{k},t_{ij}] = [y_{k},t_{ij}]=0$
si $i,j,k$ sont distincts. En particulier, l'alg\`ebre de Lie $\bar\t_{1,2}$ s'identifie \`a 
l'alg\`ebre de Lie $\f_{2}$ librement engendr\'ee par les deux g\'en\'erateurs
$x:= x_{1}$ et $y:= y_{1}$.

Soit $\HH := \{\tau\in\CC |\Im(\tau)>0\}$ le demi-plan de Poincar\'e. On note $(z,\tau)\mapsto\theta_\tau(z)$
la fonction sur $\mathbb{C}\times\HH$ donn\'ee par\footnote{On note $\on{i}:= \sqrt{-1}$.}  
$$
\theta_\tau(z):={{e^{\pi\i z}-e^{-\pi\i z}}\over{2\pi\i}}\prod_{n>0}{{(1-e^{2\pi\i(z+n\tau)})(1-e^{2\pi\i(-z+n\tau)})}
\over{(1-e^{2\pi\i n\tau})^2}}. 
$$
On a $\theta_{\tau}(z+1) = -\theta_{\tau}(z)=\theta_\tau(-z)$, $\theta_{\tau}(z+\tau) = 
-e^{-\on{i}\pi\tau}e^{-2\pi\on{i}z}\theta_{\tau}(z)$, 
${\partial\over{\partial z}}\theta_{\tau}(z)_{|z=0}=1$, et $(\theta_{\tau}(-))^{-1}(0) = \ZZ+\tau\ZZ$. 
La fonction $\theta_{\tau}(-)$ est reli\'ee \`a la fonction th\^eta de Jacobi 
donn\'ee par 
$$
\vartheta_1(z,\tau)=-\sum_{n\in{\mathbb{Z}}+{1\over 2}}e^{\pi\i n^2+2\pi\i n(z+{1\over 2})}
$$ 
par l'identit\'e
$\vartheta_1(z,\tau)=2\pi\eta(\tau)^3\theta_\tau(z)$, o\`u $\eta(\tau)=q^{{1\over 24}}\prod_{n>0}(1-q^n)$ et $q=e^{2\pi\i\tau}$.

On d\'efinit $A(\tau),B(\tau)$ comme les holonomies r\'egularis\'ees de l'\'equation diff\'erentielle
 \begin{equation} \label{eq:dept}
X'(z) =  -{{\theta_{\tau}(z+\on{ad}x)\on{ad}x}\over{\theta_{\tau}(z)
 \theta_{\tau}(\on{ad}x)}}(y)\cdot X(z), 
 \end{equation}
\`a valeurs dans le groupe $\on{exp}(\hat{\bar\t}_{1,2})$
le long les chemins $[0,1]$ et $[0,\tau]$ (l'alg\`ebre de Lie est compl\'et\'ee
pour le degr\'e en $x,y$). Dans cette \'equation, on donne \`a l'expression $f(z,{\mathrm{ad}}x)(y)$
la valeur $\sum_{n\geq 0}g_n(z)({\mathrm{ad}}x)^n(y)$, o\`u $g_n(z):={1\over{n!}}({\partial\over{\partial t}})^n f(z,t)_{|t=0}$, 
et $f(z,t)$ est une fonction analytique en deux variables, r\'eguli\`ere en $t=0$ (ici $f(z,t)=-{{\theta_\tau(z+t)}\over{\theta_\tau(z)}}
{t\over{\theta_\tau(t)}}$, qui vaut $-1$ en $t=0$). 
Cette \'equation  
admet une solution $X(z)$ d\'efinie sur $\{a+b\tau|(a,b)\in]0,1[^2\}$
telle que $X(z)\simeq (-2\pi\on{i}z)^{-[x,y]}$ en $z\to 0$. Alors 
$$
A(\tau):= X(z)^{-1}X(z+1), \quad  B(\tau):= 
X(z)^{-1}e^{2\pi\on{i}x}X(z+\tau). 
$$
Ce sont des \'el\'ements du groupe $\on{exp}(\hat{\bar\t}_{1,2})$. 

\subsection{Propri\'et\'es de $A(\tau),B(\tau)$}

\subsubsection{Propri\'et\'es modulaires} 

On a d'apr\`es \cite{En2}, Proposition 66 
\begin{equation} \label{id:mod:A:B}
A({{-1}\over{\tau}}) = \on{Ad}\big({{-1}\over\tau}\big)^{-t}
\circ \alpha_{\tau}(B(\tau)^{-1}), \quad 
B({{-1}\over{\tau}}) = \on{Ad}\big({{-1}\over{\tau}}\big)^{-1}
\circ\alpha_{\tau}(BAB^{-1}(\tau)) ,  
\end{equation}
o\`u $\alpha_{\tau}\in\on{Aut}(\on{exp}(\hat\f_{2}))$ est donn\'e par
$x\mapsto -\tau x$, $y\mapsto - 2\pi\on{i}x - \tau^{-1}y$, 
$$t:= -[x,y]$$
et $(-\tau^{-1})^{-t}:= \on{exp}(-\on{log}(-\tau^{-1})t)$, la d\'etermination 
du logarithme \'etant de partie imaginaire comprise entre $0$ et $\pi$. 

Les identit\'es $(-I_{2})(A(\tau)) = A(\tau)^{2,1}$, $(-I_{2})(B(\tau)) = B(\tau)^{2,1}$, dans lesquelles
$(-I_{2})$ est l'automorphisme de $\on{exp}(\hat\f_{2})$ induit par $x\mapsto -x$, $y\mapsto -y$, 
et les identit\'es (\ref{A:A:B:B}) (voir section \ref{rels:alg}) impliquent 
$$
(-I_{2})(A(\tau)) = e^{-\on{i}\pi t}
 A(\tau)^{-1} e^{-\on{i}\pi t}, \quad 
 (-I_{2})(B(\tau)) = e^{\on{i}\pi t}B(\tau)^{-1}e^{\on{i}\pi t}.  
$$ 
Compte tenu des ces identit\'es et de (\ref{comm:A:B}), les deux parties de (\ref{id:mod:A:B})
sont \'equivalentes apr\`es application de $\tau\mapsto -\tau^{-1}$ \`a l'une d'elles. 

\subsubsection{Relations alg\'ebriques} \label{rels:alg}

Soit $\Phi(a,b)$ l'associateur KZ, 
d\'efini par $\Phi(a,b) = Y_{1}^{-1}Y_{0}$, o\`u $Y_{i}$ sont les
solutions de $Y'(z) = (a/z+b/(z-1))Y(z)$ sur $]0,1[$ telles que 
$Y_{0}(z)\simeq z^{a}$ en $z\to 0$, $Y_{1}(z)\simeq (1-z)^{b}$ en $z\to 1$, 
et o\`u $a,b$ sont des variables formelles non-commutatives. 

On pose 
$$
\alpha_{+}:= e^{\on{i}\pi(t_{12}+t_{13})}A(\tau)^{1,23}\Phi(t_{12},t_{23}), 
\quad 
\alpha_{-}:= e^{-\on{i}\pi(t_{12}+t_{13})}B(\tau)^{1,23}\Phi(t_{12},t_{23}), 
$$
o\`u $a\mapsto a^{1,23}$ est le morphisme d'alg\`ebres de Lie 
$\bar\t_{1,2}\to\bar\t_{1,3}$ tel que $x_{1}\mapsto x_{1}$, 
$x_{2}\mapsto x_{2}+x_{3}$, $y_{1}\mapsto y_{1}$, $y_{2}\mapsto y_{2}+y_{3}$. 

La premi\`ere famille de relations satisfaites par $(A(\tau),B(\tau))$ est 
\begin{equation} \label{first:rel}
\alpha_{\pm}^{3,1,2}\alpha_{\pm}^{2,3,1}\alpha_{\pm} = 1 \quad \text{(dans }
\on{exp}(\hat{\bar\t}_{1,3})), 
\end{equation}
o\`u $a\mapsto a^{2,3,1}$ est l'automorphisme de $\bar\t_{1,3}$  
tel que $x_{i}\mapsto x_{i+1\text{ mod }3}$, $y_i\mapsto y_{i+1\text{ mod }3}$, 
et $a\mapsto a^{3,1,2}$ est le carr\'e de cet automorphisme. 

Le couple $(A(\tau),B(\tau))$ satisfait d'autre part la relation
\begin{equation} \label{second:rel}
(\Phi(t_{12},t_{23})^{-1} *  B(\tau)^{1,23},
(e^{-\on{i}\pi t_{12}}\Phi(t_{21},t_{13})) * (A(\tau)^{2,13})^{-1}) 
= e^{2\pi\on{i}t_{12}} \quad \text{(dans }
\on{exp}(\hat{\bar\t}_{1,3})), 
\end{equation}
o\`u $(x,y):= xyx^{-1}y^{-1}$, $x*y:= xyx^{-1}$, $t_{12}:= [x_{1},y_{2}]$ et 
$a\mapsto a^{2,13}$ est le morphisme $\bar\t_{1,2}\to\bar\t_{1,3}$
donn\'e par $x_{1}\mapsto x_{2}$, $x_{2}\mapsto x_{1}+x_{3}$, 
$y_{1}\mapsto y_{2}$, $y_{2}\mapsto y_{1}+y_{3}$.

Les relations (\ref{first:rel}) impliquent alors 
\begin{equation} \label{A:A:B:B}
e^{\on{i}\pi t}A(\tau)e^{\on{i}\pi t}A(\tau)^{2,1} = 
e^{-\on{i}\pi t}B(\tau)e^{-\on{i}\pi t}B(\tau)^{2,1} = 1
\quad \text{(dans }\on{exp}(\hat{\bar\t}_{1,2})), 
\end{equation}
o\`u $a\mapsto a^{2,1}$ est l'automorphisme involutif de $\overline\t_{1,2}$
donn\'e par $x_{1}\leftrightarrow x_{2}$, $y_{1}\leftrightarrow y_{2}$, 
et la relation (\ref{second:rel}) implique 
\begin{equation} \label{comm:A:B}
(A(\tau),B(\tau)) = e^{-2\pi\on{i}t} \quad \text{(dans }
\on{exp}(\hat{\bar\t}_{1,2})). 
\end{equation}

\subsubsection{Equations diff\'erentielles}\label{sect:ED}

Pour chaque $n\geq -1$, il existe une unique d\'erivation de $\f_{2}$, homog\`ene 
pour le bidegr\'e en $(x,y)$ et telle que $\delta_{2n}(x) = \on{ad}(x)^{2n+2}(y)
=:[x^{2n+2}y]$ et $\delta_{2n}([x,y])=0$. Les fonctions $A(\tau),B(\tau)$ satisfont alors les \'equations
diff\'erentielles
\begin{equation} \label{ED:A:B}
2\pi\on{i}\partial_{\tau}A(\tau) = 
-(\sum_{n\geq -1}(2n+1)G_{2n+2}(\tau)\delta_{2n})(A(\tau)), \quad 
2\pi\on{i}\partial_{\tau}B(\tau) = -(\sum_{n\geq -1}
(2n+1)G_{2n+2}(\tau)\delta_{2n})(B(\tau)). 
\end{equation}
(cf. \cite{En2}, Proposition 67), o\`u les s\'eries d'Eisenstein sont d\'efinies par 
$$
G_{k}(\tau) = \sum_{a\in(\ZZ+\tau\ZZ)-\{0\}} {1\over {a^{k}}}
\quad \text{si\ } k\text{\ est\ pair\ }\geq 4, 
\quad G_{2}(\tau) = \sum_{m\in\ZZ}(\sum_{n}{}'{1\over{(n+m\tau)^{2}}}), 
\quad G_{0}(\tau) := -1, 
$$
et $\sum'$ signifie $\sum_{n\in\ZZ}$ si $m\neq 0$ et $\sum_{n\in\ZZ - \{0\}}$
si $m=0$. 

\subsubsection{Comportement \`a l'infini} 

On a 
\begin{equation} \label{comp:A}
A(\tau) = \Phi(\tilde y,t)e^{2\pi\on{i}\tilde y} \Phi(\tilde y,t)^{-1} 
+ O(e^{2\pi\on{i}\tau}), 
\end{equation}
\begin{equation} \label{comp:B}
B(\tau)=e^{\on{i}\pi t}\Phi(-\tilde y-t,t)e^{2\pi\on{i}x}e^{2\pi\on{i}
\tilde y\tau} \Phi(\tilde y,t)^{-1} + O(e^{2\pi\on{i}(1-\epsilon)\tau})
\end{equation}
pour tout $\epsilon>0$, lorsque $\tau\to\on{i}\infty$ (\cite{CEE}, 
d\'emonstration de Proposition 4.7 puis Lemma 4.14), o\`u 
$$
\tilde y:= -{{\on{ad}x}\over{e^{2\pi\on{i}\on{ad}x}-1}}(y) 
$$  
et on rappelle que $t = -[x,y]$ et 
$\Phi(a,b)$ est l'associateur KZ d\'efini en section \ref{rels:alg}. 

\section{Les analogues elliptiques des nombres multiz\'etas} \label{sec:dev:A}

Dans cette section, nous calculons l'holonomie r\'egularis\'ee de certaines \'equations diff\'erentiel\-les 
(section \ref{subsect:IIreg}). Nous exprimons cette holonomie en termes d'int\'egrales it\'er\'ees r\'egularis\'ees 
(section \ref{subsect:IIreg:holon}).  Nous utilisons ces r\'egulari\-sa\-tions en section 
\ref{subsect:multizeta} pour d\'efinir les fonctions du param\`etre elliptique, analogues des nombres multiz\'etas. Nous 
montrons que les fonctions $A(\tau),B(\tau)$ peuvent s'interpr\'eter comme 
des s\'eries g\'en\'eratrices pour ces fonctions (section \ref{sect:corr}).
Les propri\'et\'es de $A(\tau),B(\tau)$ peuvent donc se traduire en termes
fonctionnels : c'est ce qui est fait explicitement en section \ref{subsect:proprietes} 
pour certaines de ces propri\'et\'es.  

\subsection{Holonomies r\'egularis\'ees} \label{subsect:IIreg}

Soit $I$ un ensemble fini contenant les \'el\'ements $0,1$. Soit $\Omega^1:=\Omega^1(]0,1[,{{\mathbb{C}}})$
l'espace des formes diff\'erentielles sur $]0,1[$ et soit $i\mapsto\ul\omega_i$ une application 
$I\to\Omega^1$, telle que $\ul\omega_0=d{\mathrm{log}}(z)$, $\ul\omega_1=d{\mathrm{log}}(1-z)$, 
et pour $i\neq 0,1$, la forme $\ul\omega_i$ est r\'eguli\`ere en $0$ et $1$. 

Soit $V:=\oplus_{i\in I}{\mathbb C}a_i$ l'espace vectoriel engendr\'e par $I$. On pose 
$$
\omega:=\sum_{i\in I}\ul\omega_i\cdot a_i\in\Omega^1(]0,1[,V). 
$$
On note $T(V)$ l'alg\`ebre tensorielle de $V$, munie du produit de concat\'enation, et $\hat T(V)$ sa
compl\'etion pour le degr\'e pour lequel $V$ est de degr\'e 1. L'\'equation diff\'erentielle $df=\omega f$
d'inconnue une fonction $f:]0,1[\to\hat T(V)$ admet deux solutions $f_0$ et $f_1$, telles que 
$f_0(z)\sim z^{a_0}$ pour $z\to 0$ et $f_1(z)\sim (1-z)^{a_1}$ pour $z\to 1$. On d\'efinit alors 
l'holonomie r\'egularis\'ee de cette \'equation diff\'erentielle comme 
\begin{equation}\label{def:Hol}
Hol([0,1],\omega):=f_1^{-1}f_0\in\hat T(V).
\end{equation}

D'apr\`es \cite{En1}, appendice, on a 
\begin{align*}
& Hol([0,1],\omega)\\ & =1+
\sum_{n\geq 1}\sum_{(i_1,\ldots,i_n)\in I^n, \atop i_1\neq 0,i_n\neq 1}
\Big(\int_{[0,1]}\ul\omega_{i_1}\bullet\cdots\bullet\ul\omega_{i_n}\Big)
\sum_{A\subset\{a|i_a=0\},\atop B\subset\{b|i_b=1\}}
(-a_1)^{|B|}(a_{i_n}\cdots a_{i_1})_{A,B}(-a_0)^{|A|}, 
\end{align*}
o\`u $(a_{i_n}\cdots a_{i_1})_{A,B}$ est le produit des $a_{i_k}$, dans lequel $k$ parcourt de facon 
d\'ecroissante l'ensemble $[1,n]-(A\cup B)$, et o\`u on pose 
$\int_{[0,1]}\alpha_{1}\bullet\cdots\bullet\alpha_{n}
:= \int_{\Delta_{n}}\alpha_{1}\otimes\cdots\otimes\alpha_{n}$,
o\`u $\Delta_{n}$ est le simplexe $\{(t_{1},\ldots,t_{n})\in{\mathbb{R}}^{n}|
0\leq t_{1}\leq\ldots\leq t_{n}\leq 1\}$ pour $\alpha_1,\ldots,\alpha_n\in\Omega^1$. 
 
Soit $Hol_n([0,1],\omega)$ la partie de degr\'e $n$ de cette somme. Alors 
$$
Hol([0,1],\omega)=\sum_{n\geq 0}Hol_n([0,1],\omega).
$$ 
Soit $x\mapsto x^{op}$ l'involution de $T(V)$ donn\'ee par $(v_1\otimes\cdots\otimes v_k)^{op} 
:=v_k\otimes\cdots\otimes v_1$. Soit $m:{\mathbb C}[a_0]\otimes T(V)\otimes {\mathbb C}[a_1]\to 
T(V)$ l'application compos\'ee ${\mathbb C}[a_0]\otimes T(V)\otimes {\mathbb C}[a_1]\hookrightarrow T(V)^{\otimes 
3}\to T(V)$, o\`u la premi\`ere application est un produit tensoriel d'injections canoniques et de l'identit\'e, et o\`u la 
deuxi\`eme application est la multiplication de $T(V)$. Soit $Op$ l'endomorphisme de $T(V)$ donn\'e par 
\begin{equation}\label{def:Op}
Op=(T(V)\stackrel{mor}{\to}{\mathbb C}[a_0]\otimes T(V)\otimes 
{\mathbb C}[a_1]\stackrel{m}{\to}T(V)),
\end{equation} 
o\`u $mor$ est le morphisme d'alg\`ebres induit par $a_0\mapsto a_0^{(2)}- a_0^{(1)}$, $a_1\mapsto a_1^{(2)}-a_1^{(3)}$, 
$a_i\mapsto a_i^{(2)}$, $i\neq 0$  (on note $a_0^{(1)} :=a_0\otimes 1\otimes 1$, $a_i^{(2)} :=1\otimes a_i\otimes 1$, $a_1^{(3)} :=1\otimes 
1\otimes a_1$). 

Posons $\omega_0:=\ul{\omega}_0\otimes a_0$, $\omega_1:=\ul{\omega}_1\otimes a_1$ ; ce sont des \'el\'ements 
de $\Omega^1\otimes V$. Alors 
$$
Hol_0([0,1],\omega)=1, \quad Hol_1([0,1],\omega)=(\int_{[0,1]}\otimes\mathrm{id})(\omega-\omega_0-\omega_1), 
$$
et si $n\geq 2$, 
$$
Hol_n([0,1],\omega)^{op}=(\int_{\Delta_n}\otimes Op)((\omega-\omega_0)\circ\omega^{\circ n-2}\circ(\omega-\omega_1))
$$
o\`u $\circ$ est le produit de l'alg\`ebre tensorielle $T(\Omega^1\otimes V)$. 

Notons $\omega\mapsto\omega^{01},\omega^{02},\omega^{03}$ les applications naturelles de
$\Omega^1\otimes {\mathbb C}a_0$, $\Omega^1\otimes V$, $\Omega^1\otimes {\mathbb C}a_1$ vers 
$T(\Omega^1)\otimes{\mathbb C}[a_0]\otimes T(V)\otimes{\mathbb C}[a_1]$. Alors (\ref{def:Op}) implique que pour $n\geq 2$, 
\begin{equation}\label{eq:Hol}
Hol_n([0,1],\omega)^{op}=(\int_{\Delta_n}\otimes m)\Big(((\omega-\omega_0)^{02}-\omega_1^{03})
\circ(\omega^{02}-\omega_0^{01}-\omega_1^{03})^{\circ n-2}\circ((\omega-\omega_1)^{02}-\omega_0^{01})\Big), 
\end{equation}
o\`u le terme entre crochets est un \' el\' ement de $T(\Omega^1)\otimes{\mathbb C}[a_0]\otimes T(V)\otimes{\mathbb C}[a_1]$. 

On montre facilement l'\'enonc\'e suivant : 

\begin{lemma}
Soit $\Omega,V$ des espaces vectoriels, soit $A$ une alg\`ebre commutative. On note 
$\omega\mapsto \omega^{12}$, $\alpha\mapsto\alpha^{13}$ les applications naturelles 
de $\Omega\otimes V$ et $\Omega\otimes A$ vers $T(\Omega)\otimes T(V)\otimes A$. 
On note $\sh$ le produit de battage de $T(\Omega)$, donn\'e par 
$(x_1\otimes\cdots\otimes x_k)(x_{k+1}\otimes \cdots\otimes x_{k+l})=\sum x_{f(1)}\otimes\cdots\otimes x_{f(k+l)}$, 
o\`u la somme parcourt les permutations de $[1,k+l]$ telles que si $f(i)<f(j)\leq k$ ou $l+1\leq f(i)<f(j)$, alors 
$i<j$. On note \'egalement $\sh:(T(\Omega)\otimes T(V))\otimes (T(\Omega)\otimes A)\to 
T(\Omega)\otimes T(V)\otimes A$ le produit tensoriel du produit $\sh:T(\Omega)^{\otimes 2}\to T(\Omega)$
avec les endomorphismes identit\'e de $T(V)$ et de $A$.  

Si $\omega\in \Omega\otimes V$, $\alpha\in\Omega\otimes A$, alors 
$$
(\omega^{12}+\alpha^{13})^{\circ n}=\sum_{k=0}^n \omega^{\circ n-k}\sh\alpha^{\circ k}, 
$$
o\`u les puissances sont calcul\'ees dans les produits tensoriels $T(\Omega)\otimes T(V)\otimes A$, $T(\Omega)\otimes T(V)$ et $T(\Omega)\otimes A$, o\`u $T(\Omega)$ est muni du produit de battage et $T(V)$ du produit de concat\'enation.  
\end{lemma}
En appliquant cette \'egalit\'e au membre de droite de (\ref{eq:Hol}), on trouve 
\begin{align}\label{eq:Hol:2}
Hol_n([0,1],\omega)^{op}=\sum_{k,l\geq 0\atop k+l\leq n-2} & (-a_0)^k\Big(\int_{\Delta_n}
(\omega-\omega_0)\circ(\ul\omega_0^{\circ k}\sh\ul\omega_1^{\circ l}\sh\omega^{\circ n-2-k-l})\circ(\omega-\omega_1)\Big)(-a_1)^l \nonumber
\\ & + (-a_0)^k\Big(\int_{\Delta_n}
\ul\omega_1\circ(\ul\omega_0^{\circ k}\sh\ul\omega_1^{\circ l}\sh\omega^{\circ n-2-k-l})\circ(\omega-\omega_1)\Big)(-a_1)^{l+1} \nonumber
\\ &+ (-a_0)^{k+1}\Big(\int_{\Delta_n}
(\omega-\omega_0)\circ(\ul\omega_0^{\circ k}\sh\ul\omega_1^{\circ l}\sh\omega^{\circ n-2-k-l})\circ\ul\omega_0)\Big)(-a_1)^l \nonumber
\\ & + (-a_0)^{k+1}\Big(\int_{\Delta_n}
\ul\omega_1\circ(\ul\omega_0^{\circ k}\sh\ul\omega_1^{\circ l}\sh\omega^{\circ n-2-k-l})
\circ\ul\omega_0)\Big)(-a_1)^{l+1} 
\end{align}
pour $k\geq 2$, expression dans lequelle $\circ$ est le produit de concat\'enation dans $T(\Omega^1)$ ou 
$T(\Omega^1)\otimes T(V)$, l'espace $T(\Omega^1)$ est consid\'er\'e comme un sous-espace de $T(\Omega^1)\otimes T(V)$
par tensorisation avec 1, le symbole $\sh$ d\'esigne le produit sur $T(\Omega^1)\otimes T(V)$, produit tensoriel du produit de
battage et du produit de concat\'enation, et $\int_{\Delta_n}:T(\Omega^1)\otimes T(V)\to T(V)$ d\'esigne le produit tensoriel de 
$\int_{\Delta_n}:T(\Omega^1)\to{\mathbb C}$ avec l'identit\'e de $T(V)$. 

En utilisant l'identit\'e $(a\circ A)\sh(b\circ B)=a\circ (A\sh(b\circ B))
+b\circ(B\sh(a\circ A))$ dans $T(\Omega^1)$, o\`u $a,b\in\Omega^1$, $A,B\in T(\Omega^1)$, 
on simplifie ainsi cette expression 
\begin{equation}\label{simplified}
Hol_n([0,1],\omega)^{op}=\int_{\Delta_n}
\sum_{k,l\geq 0\atop k+l\leq n-2} (-a_0)^k\big(\ul\omega_0^{\circ k}
\sh\ul\omega_1^{\circ l}\sh\omega^{\circ n-2-k-l}\big)(-a_1)^l, 
\end{equation}
dans laquelle les signes somme et int\'egrale ne peuvent \^etre invers\'es, les termes individuels de la somme n'\'etant pas int\'egrables. 

\begin{remark}
On peut r\'eduire la dimension du simplexe d'int\'egration dans les formules (\ref{eq:Hol:2}), (\ref{simplified}) en utilisant l'identit\'e
$$
\int_{\Delta_n} \alpha\circ ((df)^{\circ a}\sh(dg)^{\circ b}\sh\beta)\circ\gamma
=\sum_{a',a''|a'+a''=a}\sum_{b',b''|b'+b''=b}\int_{\Delta_{n-a-b}}
{(-f)^{a''}\over{a''!}}{(-g)^{b''}\over{b''!}}\alpha\circ\beta\circ\gamma{f^{a'}\over{a'!}} {g^{b''}\over{b''!}}
$$
dans laquelle $\alpha,\beta,\gamma$ sont dans $T(\Omega^1)$ et $f,g$ sont des fonctions sur $[0,1]$. 
\end{remark}

\subsection{Int\'egrales it\'er\'ees r\'egularis\'ees}\label{subsect:IIreg:holon}

Soit $\mathcal{A}$ une alg\`ebre, limite projective d'alg\`ebres de dimension finie. Soit 
${\mathbf{a}}_0$, ${\mathbf{a}}_1$ des \'el\'ements de  $\mathcal{A}$, et posons 
\begin{equation}\label{def:Omega}
\Omega_{{\mathbf{a}}_0,{\mathbf{a}}_1}:=\{
\mbox{\boldmath$\omega$}
\in\Omega^1(]0,1[,{\mathcal{A}})|\mbox{\boldmath$\omega$}={\mathbf{a}}_0d{\mathrm{log}}(z)+O(1)\text{\ si\ }z\to 0,\ 
\mbox{\boldmath$\omega$}={\mathbf{a}}_1d{\mathrm{log}}(1-z)+O(1)\text{\ si\ }z\to 1\}. 
\end{equation}
Pour $\mbox{\boldmath$\omega$}\in\Omega_{{\mathbf{a}}_0,{\mathbf{a}}_1}$, on pose 
$$
I^{reg}_{[0,1]}(\mbox{\boldmath$\omega$}):=\int_{[0,1]}(\mbox{\boldmath$\omega$}-{\mathbf{a}}_0d{\mathrm{log}}(z)
-{\mathbf{a}}_1d{\mathrm{log}}(1-z)).
$$ 
Pour $n\geq 2$ et $\mbox{\boldmath$\omega$}_1,\ldots,\mbox{\boldmath$\omega$}_n\in\Omega_{{\mathbf{a}}_0,{\mathbf{a}}_1}$, on pose 
$$
I^{reg}_{[0,1]}(\mbox{\boldmath$\omega$}_1,\ldots,\mbox{\boldmath$\omega$}_n):=
\sum_{a,b,\epsilon,\eta|\atop{a,b\geq 0, a+b\leq n-2,\atop{\epsilon,\eta\in\{0,1\}}}} {\mathbf{a}}_0^{a+\epsilon}\cdot
I_{a,b}^{\epsilon\eta} \cdot {\mathbf{a}}_1^{b+\eta},
$$
le produit \'etant calcul\'e dans ${\mathcal{A}}^{op}$ (l'alg\`ebre oppos\'ee \`a ${\mathcal{A}}$), o\`u 
$$
I^{00}_{a,b}=\int_{\Delta_n}(\mbox{\boldmath$\omega$}_{a+1}-{\mathbf{a}}_0d{\mathrm{log}}(z))
\circ\Big((-d{\mathrm{log}}(z))^{\circ a}\sh(-d{\mathrm{log}}(1-z))^{\circ b}\sh
(\mbox{\boldmath$\omega$}_{a+2}\circ\cdots\circ\mbox{\boldmath$\omega$}_{n-b-1})\Big)\circ
(\mbox{\boldmath$\omega$}_{n-b}-{\mathbf{a}}_1d{\mathrm{log}}(1-z)), 
$$
$$
I^{01}_{a,b}=\int_{\Delta_n}(-d{\mathrm{log}}(1-z))
\circ\Big((-d{\mathrm{log}}(z))^{\circ a}\sh(-d{\mathrm{log}}(1-z))^{\circ b}\sh
(\mbox{\boldmath$\omega$}_{a+1}\circ\cdots\circ\mbox{\boldmath$\omega$}_{n-b-2})\Big)\circ
(\mbox{\boldmath$\omega$}_{n-b-1}-{\mathbf{a}}_1d{\mathrm{log}}(1-z)), 
$$
$$
I^{10}_{a,b}=\int_{\Delta_n}(\mbox{\boldmath$\omega$}_{a+2}-{\mathbf{a}}_0d{\mathrm{log}}(z))
\circ\Big((-d{\mathrm{log}}(z))^{\circ a}\sh(-d{\mathrm{log}}(1-z))^{\circ b}\sh
(\mbox{\boldmath$\omega$}_{a+3}\circ\cdots\circ\mbox{\boldmath$\omega$}_{n-b})\Big)\circ
(-d{\mathrm{log}}(z)), 
$$
$$
I^{11}_{a,b}=\int_{\Delta_n}(-d{\mathrm{log}}(1-z))
\circ\Big((-d{\mathrm{log}}(z))^{\circ a}\sh(-d{\mathrm{log}}(1-z))^{\circ b}\sh
(\mbox{\boldmath$\omega$}_{a+2}\circ\cdots\circ\mbox{\boldmath$\omega$}_{n-b-1})\Big)\circ
(-d{\mathrm{log}}(z)), 
$$
o\`u : $\circ$ est le produit dans $T(\Omega^1)\otimes{\mathcal{A}}^{op}$, les formes
$d{\mathrm{log}}(z)$, $d{\mathrm{log}}(1-z)$ de $\Omega^1$ sont identifi\'ees \`a des \'el\'ements de 
$T(\Omega^1)\otimes{\mathcal{A}}^{op}$ par tensorisation avec $1$, on note l'application 
$\int_{\Delta_k}\otimes{\mathrm{id}}:(\Omega^1)^{\otimes k}\otimes{\mathcal{A}}^{op}\to{\mathcal{A}}^{op}$
simplement $\int_{\Delta_k}$.

On a alors : 

\begin{proposition} Soit $\mathcal{A}$ une alg\`ebre, limite projective d'alg\`ebres de dimension finie. 
Soit ${\mathbf{a}}_0,{\mathbf{a}}_1$ des \'el\'ements de $\mathcal{A}$ et soit $\mbox{\boldmath$\omega$}$ un \'el\'ement de 
$\Omega_{{\mathbf{a}}_0,{\mathbf{a}}_1}$. L'holonomie r\'egularis\'ee $Hol([0,1],\mbox{\boldmath$\omega$})$ de 
l'\'equation diff\'erentielle $df= \mbox{\boldmath$\omega$}f$, d\'efinie par 
(\ref{def:Hol}), est donn\'ee par 
\begin{equation}\label{form:Hol}
Hol([0,1],\mbox{\boldmath$\omega$})=1+\sum_{n\geq 1}I_{[0,1]}^{reg}(\underbrace{\mbox{\boldmath$\omega$},\ldots,
\mbox{\boldmath$\omega$}}_n). 
\end{equation}   
\end{proposition}

Soit $M$ une vari\'et\'e lisse, $U\subset M$ un ouvert, et $\gamma:[0,1]\to M$ un chemin tel que $\gamma(]0,1[)\subset U$. 
Pour ${\mathcal{A}},{\mathbf{a}}_0,{\mathbf{a}}_1$ comme ci-dessus, on pose 
\begin{equation}\label{def:Omega:chemin}
\Omega_{{\mathbf{a}}_0,{\mathbf{a}}_1}(\gamma):=\{\mbox{\boldmath$\omega$}\in
\Omega^1(U,{\mathcal{A}})|\gamma^*(\mbox{\boldmath$\omega$})\in\Omega_{{\mathbf{a}}_0,{\mathbf{a}}_1}\} ;
\end{equation} 
pour $\mbox{\boldmath$\omega$}_1,\ldots,\mbox{\boldmath$\omega$}_n\in\Omega_{{\mathbf{a}}_0,{\mathbf{a}}_1}(\gamma)$, on pose $I^{reg}_\gamma(\mbox{\boldmath$\omega$}_1,\ldots,\mbox{\boldmath$\omega$}_n):=I^{reg}_{[0,1]}
(\gamma^*(\mbox{\boldmath$\omega$}_1),\ldots,\gamma^*(\mbox{\boldmath$\omega$}_n))$. Enfin, pour 
$\mbox{\boldmath$\omega$}\in\Omega_{{\mathbf{a}}_0,{\mathbf{a}}_1}(\gamma)$, on pose $Hol(\gamma,\mbox{\boldmath$\omega$}):=Hol([0,1],\gamma^*(\mbox{\boldmath$\omega$}))$. On a alors 
$$
Hol(\gamma,\mbox{\boldmath$\omega$})=1+\sum_{n\geq 1}I_\gamma^{reg}(\underbrace{\mbox{\boldmath$\omega$},\ldots,
\mbox{\boldmath$\omega$}}_n).  
$$

\begin{example} For $n=2$, 
\begin{align*}
I^{reg}_{[0,1]}(\mbox{\boldmath$\omega$}_1,\mbox{\boldmath$\omega$}_2)=
\int_{\Delta_2} & (\mbox{\boldmath$\omega$}_1-{\mathbf{a}}_0d{\mathrm{log}}(z))\circ(\mbox{\boldmath$\omega$}_2-{\mathbf{a}}_1d{\mathrm{log}}(1-z)) 
+ (-d{\mathrm{log}}(1-z))\circ(\mbox{\boldmath$\omega$}_2-{\mathbf{a}}_1d{\mathrm{log}}(z))\cdot{\mathbf{a}}_1
\\ & +{\mathbf{a}}_0\cdot(\mbox{\boldmath$\omega$}_1-{\mathbf{a}}_0d{\mathrm{log}}(z))\circ(-d{\mathrm{log}}(z))
+{\mathbf{a}}_0\cdot(-d{\mathrm{log}}(1-z))\circ(-d{\mathrm{log}}(z))\cdot{\mathbf{a}}_1. 
\end{align*}
\end{example}

\subsection{Les analogues elliptiques des nombres multiz\'etas}\label{subsect:multizeta}

Fixons $\tau\in\HH$. Pour $x\in\CC$, on pose 
$$
\sigma^{\tau}_{x}(z):= {\theta_\tau(z+x)\over{\theta_\tau(z)\theta_\tau(x)}}. 
$$
Consid\'erant $x$ comme une variable formelle proche de $0$, 
on voit $\sigma_{x}^\tau$ comme un \'el\'ement de  $x^{-1}\on{Mer}(\CC)[[x]]$, 
o\`u $\on{Mer}(\CC) = \{$fonctions m\'eromorphes d\'efinies sur $\CC\}$. 
Plus pr\'ecis\'ement : 

\begin{proposition}
$\sigma_{x}^{\tau}$ admet le d\'eveloppement 
$$
\sigma_{x}^{\tau}(z) = {1\over x} + \sum_{n\geq 0}k_{n}^{\tau}(z)x^{n}, 
$$
avec $k^{\tau}_{0}(z) = (\theta'_{\tau}/\theta_{\tau})(z)$ et $k_{n}^{\tau}$ 
finie en $0$ et $1$ si $n>0$. 
\end{proposition}

{\em D\'emonstration.} Le param\`etre $\tau$ \'etant fix\'e, on consid\`ere $\theta_\tau(-)$ comme 
une fonction de la variable $z$.  
On a $x\sigma^{\tau}_{x}(z)_{|x=0}=1$. De plus, 
$$
(\sigma^{\tau}_{x}(z) - {1\over x})_{|x=0} = {1\over x}
({\theta_\tau(z+x)\over{\theta_\tau(z)}}{x\over{\theta_\tau(x)}}-1)_{|x=0} 
= {\theta'_\tau\over\theta_\tau}(z). $$
Enfin, le d\'eveloppement de $\sigma_{x}^{\tau}(z)$ en $z=0$ est, compte tenu de 
l'imparit\'e de $\theta_\tau$
$$
\sigma_{x}^{\tau}(z) = {\theta_\tau(x+z)\over{\theta_\tau(x)}}{1\over{\theta_\tau(z)}}
= (1+z{\theta'_\tau\over\theta_\tau}(x)+O(z^{2}))({1\over z}+O(z)) = 
{1\over z}+{\theta'_\tau\over\theta_\tau}(x)+O(z). 
$$ 
Donc 
$$
\sigma_{x}^{\tau}(z) - ({1\over x}+k_{0}^{\tau}(z)) 
= ({\theta'_\tau\over\theta_\tau}(x)-{1\over x})
+ ({1\over z} - {\theta'_\tau\over\theta_\tau}(z))+O(z)$$
qui est fini en $z=0$. Il s'ensuit que les $k_{n}^{\tau}$ sont finis en $0$. Par sym\'etrie, 
ils sont \'egalement finis en $1$. \hfill \qed\medskip 

Posons ${\mathcal A}:={\mathbb{C}}$, ${\mathbf{a}}_0={\mathbf{a}}_1:=1$. Alors pour $x\in{\mathbb{C}}$, 
la forme $\sigma_x^\tau(z)dz$
appartient \`a l'espace $\Omega_{1,1}$ d\'efini par (\ref{def:Omega}). 

De m\^eme, la forme 
$e^{2\pi\i{{xz}\over\tau}}\sigma_x^\tau(z)dz$ appartient \`a l'espace $\Omega_{1,1}(\underline{[0,\tau]})$ d\'efini 
par (\ref{def:Omega:chemin}), o\`u $\underline{[0,\tau]}$ est le chemin lin\'eaire $[0,1]\to[0,\tau]$ trac\'e sur $\mathbb{C}$. 

On pose alors : 

\begin{definition} On note $\underline{[0,1]}$ et $\underline{[0,\tau]}$ les chemins lin\'eaires $[0,1]\to[0,1]$ 
et $[0,1]\to[0,\tau]$ trac\' es sur $\CC$. Pour $\tau\in\HH$, on pose
\begin{equation} \label{def:I}
I_{x_{1},\ldots,x_{n}}(\tau):= 
I^{reg}_{\underline{[0,1]}}(\sigma_{x_{1}}^{\tau}dz,
\ldots,\sigma_{x_{n}}^{\tau}dz)
\end{equation}
\begin{equation} \label{def:J}
J_{x_{1},\ldots,x_{n}}(\tau):= I_{\underline{[0,\tau]}}^{reg}
(e^{2\pi\on{i}{{x_{1}z}\over\tau}}
\sigma_{x_{1}}^{\tau}(z)dz,\ldots,
e^{2\pi\on{i}{{x_{n}z}\over\tau}}\sigma_{x_{n}}^{\tau}(z)dz) ;   
\end{equation}
ce sont des s\'eries dans $(x_{1}\cdots x_{n})^{-1}\CC[[x_{1},\ldots,x_{n}]]$.   
\end{definition}
En utilisant (\ref{eq:Hol:2}), on trouve pour $n\geq 0$ 
\begin{align}\label{explicit:I}
& I_{x_1,\ldots,x_n}(\tau)=\sum_{k,l|k,l\geq 0,\atop{k+l\leq n-2}}
\\ & \nonumber \int_{\Delta_{n}}
(\sigma_{x_{k+1}}^\tau-d{\mathrm{log}}(z))\circ\Big(
(-d{\mathrm{log}}(z))^{\circ a}\sh(-d{\mathrm{log}}(1-z))^{\circ b}\sh(\sigma_{x_{k+2}}^\tau\circ\cdots\circ\sigma_{x_{n-l-1}}^\tau)\Big)\circ(\sigma_{x_{n-l}}^\tau-d{\mathrm{log}}(1-z))
\\ & \nonumber+ \int_{\Delta_{n}}(-d{\mathrm{log}}(1-z))\circ\Big(
(-d{\mathrm{log}}(z))^{\circ a}\sh(-d{\mathrm{log}}(1-z))^{\circ b}\sh(\sigma_{x_{k+1}}^\tau\circ\cdots\circ\sigma_{x_{n-l-2}}^\tau)\Big)\circ(\sigma_{x_{n-l-1}}^\tau-d{\mathrm{log}}(1-z))
\\ &  \nonumber
+ \int_{\Delta_{n}}(\sigma_{x_{k+2}}^\tau-d{\mathrm{log}}(z))\circ\Big(
(-d{\mathrm{log}}(z))^{\circ a}\sh(-d{\mathrm{log}}(1-z))^{\circ b}\sh(\sigma_{x_{k+3}}^\tau\circ\cdots\circ\sigma_{x_{n-l}}^\tau)\Big)
\circ(-d{\mathrm{log}}(z))
\\ &  \nonumber
+ \int_{\Delta_{n}}(-d{\mathrm{log}}(1-z))\circ\Big(
(-d{\mathrm{log}}(z))^{\circ a}\sh(-d{\mathrm{log}}(1-z))^{\circ b}\sh(\sigma_{x_{k+2}}^\tau\circ\cdots\circ\sigma_{x_{n-l-1}}^\tau)\Big)
\circ(-d{\mathrm{log}}(z))
\end{align}
ce qui donne 
\begin{equation} \label{condensed:I}
I_{x_1,\ldots,x_n}(\tau)=\int_{\Delta_{n}}\sum_{k,l|k,l\geq 0,\atop{k+l\leq n-2}}
(-d{\mathrm{log}}(z))^{\circ k}\sh(-d{\mathrm{log}}(1-z))^{\circ l}\circ(\sigma_{x_{k+1}}^\tau\circ\cdots\circ\sigma_{x_{n-l}}^\tau)
\end{equation}
(en notant $\sigma_x^\tau$ \`a la place de $\sigma_x^\tau(z)dz$).
En particulier, on a 
\begin{align*}
I_{x,y}(\tau)=\int_{\Delta_2}
& (\sigma_x^\tau-d{\mathrm{log}}(z))\circ(\sigma_y^\tau-d{\mathrm{log}}(1-z))-d{\mathrm{log}}(1-z)\circ(\sigma_x^\tau-d{\mathrm{log}}(1-z))
\\ & -(\sigma_y^\tau-d{\mathrm{log}}(z))\circ d{\mathrm{log}}(z)+d{\mathrm{log}}(1-z)\circ d{\mathrm{log}}(z)
\end{align*}
On a aussi 
$$
I_x(\tau)=\int_{[0,1]}(\sigma_{x}^\tau(z)dz-d{\mathrm{log}}(z)-d{\mathrm{log}}(1-z))
$$
La fonction  $J_{x_1,\ldots,x_n}(\tau)$ est donn\'ee par les formules analogues, obtenues au moyen des substitutions 
$[0,1]\to[0,\tau]$, ${\mathrm{log}}(z)\to{\mathrm{log}}(z/\tau)$, ${\mathrm{log}}(1-z)\to{\mathrm{log}}(1-{z\over\tau})$,
$\sigma_x^\tau(z)dz\to e^{2\pi\i{{xz}\over\tau}}\sigma_x^\tau(z)dz$. 

Pour $\ul{d}:= (d_{1},\ldots,d_{n})\in \ZZ_{\geq -1}^{n}$, 
on note $I_{\ul{d}}(\tau)$, $J_{\ul{d}}(\tau)$ les nombres complexes d\'efinis par 
$$
I_{\ul{x}}(\tau)=\sum_{\ul{d}\in\ZZ_{\geq -1}^{n}}I_{\ul{d}}(\tau)\ul{x}^{\ul{d}},\quad 
J_{\ul{x}}(\tau)=\sum_{\ul{d}\in\ZZ_{\geq -1}^{n}}J_{\ul{d}}(\tau)\ul{x}^{\ul{d}},
$$ 
o\`u $\ul{x}:=(x_1,\ldots,x_n)$ et $\ul{x}^{\ul{d}}:=x_1^{d_1}\cdots x_n^{d_n}$. On a en particulier 
$I_{\underbrace{\scriptstyle{-1,\ldots,-1}}_{n}}(\tau) = 1/n!$, 
$J_{\underbrace{\scriptstyle{-1,\ldots,-1}}_{n}}(\tau) = \tau^{n}/n!$.

On appelle les fonctions $I_{\ul{d}}(\tau)$, $J_{\ul{d}}(\tau)$ les 
{\it analogues elliptiques des nombres multiz\'etas} ; 
les $I_{x_{1},\ldots,x_{n}}(\tau)$, 
$J_{x_{1},\ldots,x_{n}}(\tau)$ en sont des s\'eries g\'en\'eratrices. 
Cette terminologie est justifi\'ee par les r\'esultats de la section suivante. 

\subsection{Lien avec le d\'eveloppement de $A(\tau)$, $B(\tau)$} \label{sect:corr}

Soit
$$ F:= \oplus_{n\geq 0} F_{n} := \bigoplus_{n\geq 0} (x_{1}\cdots x_{n})^{-1}
\CC[[x_{1},\ldots,x_{n}]] ;  $$
munie du produit $f*g:= h$, avec $h(x_{1},\ldots,x_{n+m}):= 
f(x_{1},\ldots,x_{n})g(x_{n+1},\ldots,x_{n+m})$ pour $f\in F_{n}$, 
$g\in F_{m}$, c'est une alg\`ebre gradu\'ee. 

L'alg\`ebre de Lie ${\mathfrak{f}}_{2}$ est librement engendr\'ee par les \'el\'ements $x,y$ et est bigradu\'ee par le degr\'e en 
ces g\'en\'erateurs. Par \'elimination de Lazard, la somme directe ${\mathfrak{f}}_{2}\ominus \CC x$ de ses composantes
de degr\'e $>0$ en $x$ est l'alg\`ebre de Lie librement engendr\'ee par les $[x^{n}y]:=({\mathrm{ad}}x)^n(y)$, $n\geq 0$. Son alg\`ebre 
enveloppante est donc l'alg\`ebre associative libre sur ces g\'en\'erateurs. On note $\mathcal{A}$ la compl\'etion de cette 
alg\`ebre enveloppante pour le degr\'e total en $x,y$. 

On en d\'eduit : 
\begin{lemma}\label{lemme:iso}
On a un isomorphisme entre $\mathcal{A}$ et la compl\' etion $\hat F:=\prod_{n\geq 0}F_n$ 
de $F$ via 
$$
[x^{d_{1}}y]\cdots [x^{d_{n}}y] \leftrightarrow 
x_{1}^{d_{1}-1}\cdots x_{n}^{d_{n}-1}\in F_{n}. 
$$  
\end{lemma}

On consid\`ere l'alg\`ebre enveloppante de ${\mathfrak{f}}_{2}\ominus \CC x$ comme gradu\'ee par le degr\'e total en 
$x,y$ ; son alg\`ebre enveloppante est donc gradu\'ee avec composantes homog\`enes de dimension finie, et $\mathcal{A}$ est donc 
une limite projective d'alg\`ebres de dimension finie. 

Posons ${\mathbf{a}}_0={\mathbf{a}}_1:=t=-[x,y]$ et 
$$
\mbox{\boldmath$\omega$}(z)dz:=-{{\theta_\tau(z+{{\mathrm{ad}}(x)})}\over{\theta_\tau(z)}}
{{{\mathrm{ad}}(x)}\over{\theta_\tau({\mathrm{ad}}(x))}}(y)dz\in\Omega^1(]0,1[,{\mathcal{A}}). 
$$
Alors $\mbox{\boldmath$\omega$}(z)dz\in\Omega_{{\mathbf{a}}_0,{\mathbf{a}}_1}$. On lui associe
$$
I(\tau):=Hol([0,1],\mbox{\boldmath$\omega$})\in{\mathcal{A}}. 
$$
On pose $\tilde{\mbox{\boldmath$\omega$}}(z)dz:=e^{2\pi\i{z\over\tau}}\mbox{\boldmath$\omega$}(z)dz$. 
Cette forme admet comme p\^ole $t$ en $0$ et $\tau$, donc $\tilde{\mbox{\boldmath$\omega$}}\in\Omega_{{\mathbf{a}}_0,{\mathbf{a}}_1}
(\underline{[0,\tau]})$, o\`u $\underline{[0,\tau]}$ est le chemin lin\'eaire $[0,1]\to[0,\tau]$. On pose alors 
$$
J(\tau):=Hol(\underline{[0,\tau]},\tilde{\mbox{\boldmath$\omega$}})\in{\mathcal{A}}. 
$$
Sous l'isomorphisme du lemme \ref{lemme:iso}, on a 
$$
{\mathcal A}\ni\mbox{\boldmath$\omega$}\leftrightarrow -\sigma^\tau_x(z)dz\in F_1, 
$$
d'o\`u
$$
{\mathcal A}\ni I^{reg}(\underbrace{\mbox{\boldmath$\omega$},\ldots,\mbox{\boldmath$\omega$}}_n)
\leftrightarrow (-1)^n I^{reg}_{[0,1]}(\sigma^\tau_{x_n},\ldots,\sigma^\tau_{x_1})=(-1)^n I_{x_n,\ldots,x_1}(\tau)\in
F_n, 
$$
d'o\`u
\begin{equation}\label{corresp:I}
{\mathcal A}\ni I(\tau)\leftrightarrow\Big((-1)^n I_{x_n,\ldots,x_1}(\tau)\Big)_{n\geq 0}\in\hat F. 
\end{equation}
On montre de m\^eme la correspondance 
\begin{equation}\label{corresp:J}
{\mathcal A}\ni J(\tau)\leftrightarrow\Big((-1)^n J_{x_n,\ldots,x_1}(\tau)\Big)_{n\geq 0}\in\hat F. 
\end{equation}
Les solutions $f_0,f_1$ de $df=\mbox{\boldmath$\omega$}f$ (voir section \ref{subsect:IIreg}) sont reli\'ees \`a $X(z)$ 
(voir section \ref{sec:def:A:B})
par $X(z)=f_0(z)(-2\pi\i)^t$, $X(z-1)=f_1(z)(2\pi\i)^t$. On en d\'eduit
\begin{equation}\label{interm:1}
{\mathrm{Ad}}((2\pi)^t)(e^{\i{\pi\over 2}t}A(\tau)e^{\i{\pi\over 2}t})=Hol([0,1],\Omega). 
\end{equation}
On a $B(\tau)=Z(z)^{-1}Z(z+\tau)$, 
o\`u $Z(z)$ est la solution dans ${\mathbf{D}}:=\{a+b\tau|(a,b)\in]0,1[^2\}$ de 
$$
dZ=({{2\pi\i x}\over{\tau}}dz+\tilde{\mbox{\boldmath$\omega$}}(z)dz)Z 
$$
telle que $Z(z)\sim(-2\pi\i z)^t$ pour $z\to 0$. Soit $e_{+}\in\on{Der}(\f_{2})$ la d\'erivation 
$(x,y)\mapsto(0,x)$. En appliquant l'automorphisme $\on{exp}({{2\pi\on{i}}\over\tau}e_{+})$ \`a l'expression reliant 
$B(\tau)$ avec $Z(z)$, on obtient : 
$$
\on{exp}({{2\pi\on{i}}\over\tau}e_{+})(B(\tau))=T(z)^{-1}T(z+\tau), 
$$
o\`u $T(z)$ est la solution dans ${\mathbf{D}}$ de $dT=\tilde{\mbox{\boldmath$\omega$}}(z)dz\cdot T$, telle que $T(z)\sim(-2\pi\i z)^t$
en $z\to 0$. Un raisonnement analogue au pr\'ec\'edent donne alors 
\begin{equation}\label{interm:2}
{\mathrm{Ad}}((2\pi)^t)(e^{-\i{\pi\over 2}t}\cdot\on{exp}({{2\pi\on{i}}\over\tau}e_{+})(B(\tau))\cdot e^{-\i{\pi\over 2}t})=Hol([0,\tau],\tilde\Omega)
\end{equation}
(\ref{corresp:I}) et (\ref{corresp:J}) d'une part, (\ref{interm:1}) et (\ref{interm:2}) d'autre part impliquent : 

\begin{proposition}
Sous l'isomorphisme du Lemme \ref{lemme:iso}, on a 
$$
{\mathrm{Ad}}((2\pi)^t)\big(e^{\i{\pi\over 2}t}A(\tau)e^{\i{\pi\over 2}t}\big)
\leftrightarrow
\Big(I_{x_n,\ldots,x_1}(\tau)\Big)_{n\geq 0}. 
$$
$$
{\mathrm{Ad}}((2\pi)^t)\big(e^{-\i{\pi\over 2}t}\cdot\on{exp}({{2\pi\on{i}}\over\tau}e_{+})(B(\tau))\cdot e^{-\i{\pi\over 2}t}\big)
\leftrightarrow
\Big((-1)^{n}J_{x_n,\ldots,x_1}(\tau)\Big)_{n\geq 0}. 
$$
\end{proposition}

\begin{remark}
Le formalisme d\'evelopp\'e par J. Ecalle utilise une variante de l'isomorphime $F\leftrightarrow\bigoplus_{n\geq 0}
{\mathbb{C}}[x_1,\ldots,x_n]$ du lemme \ref{lemme:iso}, donn\'e par $[x^{d_1}y]\cdots[x^{d_n}y]=f\leftrightarrow 
ma_f:=x_1^{d_1}\cdots x_n^{d_n}\in{\mathbb{C}}[x_1,\ldots,x_n]$.

\end{remark}

\subsection{Propri\'et\'es alg\'ebriques et modulaires} \label{subsect:proprietes}

Les multiz\'etas elliptiques sont reli\'es par l'identit\'e modulaire 
\begin{equation} \label{id:mod}
J_{x_{1},\ldots,x_{n}}(\tau) = I_{{x_{1}\over\tau},\ldots,{x_{n}\over\tau}}({{-1}\over\tau}), 
\end{equation}
qui repose sur l'identit\'e $\sigma^{-1/\tau}_{x}(z)dz = e^{2\pi\on{i}x\tau z}
\sigma^{\tau}_{\tau x}(\tau z)d(\tau z)$. 
L'identit\'e (\ref{id:mod}) traduit l'identit\'e (\ref{id:mod:A:B}) reliant
$A(\tau)$ et $B(\tau)$. 

Le caract\`ere ``de type groupe'' de $A(\tau),B(\tau)$ se traduit par les identit\'es
\begin{equation} \label{rels:Id}
I_{d_{1},\ldots,d_{n}}(\tau)
I_{d_{n+1},\ldots,d_{n+m}}(\tau) = \sum_{\sigma\in S_{n,m}}
I_{d_{\sigma(1)},\ldots,d_{\sigma(n+m)}}(\tau),
\end{equation}
o\`u $S_{n,m} = \{\sigma\in S_{n+m}|\sigma(i)<\sigma(j)$ si $i<j\leq n$ ou 
$n+1\leq i<j\}$, ainsi que les identit\'es similaires pour les $J_{\ul{d}}(\tau)$ (qui leur sont \'equivalentes 
compte tenu de (\ref{id:mod})).  

Les identit\'es (\ref{A:A:B:B}) se traduisent par
$$
\sum_{k=0}^{n} (-1)^{d_{1}+\cdots+d_{k}} I_{d_{1},\ldots,d_{k}}(\tau)
I_{d_{k+1},\ldots,d_{n}}(\tau)=0 \quad\text{si}\quad 
n\geq 1, \ d_{1},\ldots,d_{n}\geq -1$$
et les identit\'es analogues en rempla\c{c}ant chaque $I_{\ul{d}}$ par
$J_{\ul{d}}$ (qui leur sont \'egalement \'equivalentes). 

\section{Syst\`eme diff\'erentiel pour les analogues elliptiques des nombres 
multiz\'etas} \label{sect:diff}

Dans cette section, on montre que les int\'egrales it\'er\'ees r\'egularis\'ees d\'efinies dans la section \ref{subsect:IIreg:holon} sont, sous certaines 
hypoth\`eses, solutions de certains syst\`emes diff\'erentiels (section \ref{sec:3:1}). En section \ref{sec:3:2}, on applique ce r\'esultat aux 
analogues elliptiques des nombres multiz\'etas d\'efinis en section \ref{subsect:multizeta}. 

\subsection{Syst\`emes diff\'erentiels pour les int\'egrales it\'er\'ees r\'egularis\'ees}\label{sec:3:1}

Soit $I$ un intervalle de ${\mathbb R}$ contenant $0$. On note $\mathcal F$ l'ensemble des 1-formes lisses sur $I\times]0,1[$
de la forme $\omega=\omega(t,z)dz$, telles que $\omega-d{\mathrm{log}}(z)$ admet un prolongement lisse \`a $I\times[0,1[$
et $\omega-d{\mathrm{log}}(1-z)$ admet un prolongement lisse \`a $I\times]0,1]$. On note $\mathcal G$ l'ensemble des 1-formes lisses
sur $I\times[0,1]$ de la forme $g=g(t,z)dt$. On suppose donn\'es : 

$\bullet$ des \'el\'ements $\omega_i$ ($i=1,\ldots,n$) et $\psi_{i,i+1}$ ($i=1,\ldots,n-1$) de ${\mathcal F}$; 

$\bullet$ des \'el\'ements $g_i$ ($i=1,\ldots,n$) de $\mathcal G$, tels que:
\begin{itemize}
\item[(a)] pour $i=1,\ldots,n$, la forme $\omega_i(t,z)dz+g_i(t,z)dt$ est ferm\'ee, ce qui se traduit par l'identit\'e
$$
{{\partial\omega_i}\over{\partial t}}(t,z)={{\partial g_i}\over{\partial z}}(t,z)
$$
sur $I\times]0,1[$; 
\item[(b)] pour $i=1,\ldots,n-1$, on a l'identit\'e
$$
g_i(t,z)\omega_{i+1}(t,z)-g_{i+1}(t,z)\omega_i(t,z)=(g_i(t,0)-g_{i+1}(t,0))\psi_{i,i+1}(t,z)
$$
sur $I\times]0,1[$;
\item[(c)] pour $i=1,\ldots,n$, on a l'identit\'e $g_i(t,0)=g_i(t,1)$ sur $I$. 
\end{itemize}

On se place dans le cadre de la section \ref{subsect:IIreg:holon}, avec ${\mathcal{A}}={\mathbb{C}}$, ${\mathbf{a}}_0={\mathbf{a}}_1=1$, 
pour d\'efinir $\Omega_{1,1}$. Pour tout $\alpha\in{\mathcal{F}}$ et tout $t\in I$, la restriction $\alpha^t$ de $\alpha$ \`a $\{t\}\times ]0,1[$
appartient \`a $\Omega_{1,1}$. La d\'efinition de $I^{reg}_{[0,1]}$ en section \ref{subsect:IIreg:holon} permet alors de d\'efinir le nombre 
$I^{reg}_{[0,1]}(\alpha_1^t,\ldots,\alpha_n^t)$ pour tous $\alpha_1,\ldots,\alpha_n\in{\mathcal F}$ et $t\in I$. 

\begin{proposition}\label{prop:3:1}
Sous les hypoth\`eses du d\'ebut de la section \ref{sec:3:1} sur $(\omega_i)_{i=1,\ldots,n}$, $(\psi_{i,i+1})_{i=1,\ldots,n-1}$, et
$(g_i)_{i=1,\ldots,n}$, on a l'identit\'e suivante sur $I$
\begin{align}\label{eq:diff}
\nonumber {d\over{dt}}I^{reg}_{[0,1]}(\omega_1^t,\ldots,\omega_n^t)=
& -g_1^t(0)I^{reg}_{[0,1]}(\omega_2^t,\ldots,\omega_n^t)+g_n^t(0)I^{reg}_{[0,1]}(\omega_1^t,\ldots,\omega_{n-1}^t)
\\ & +\sum_{i=1}^{n-1}(g_i^t(0)-g_{i+1}^t(0))I^{reg}_{[0,1]}(\omega_1^t,\ldots,\psi_{i,i+1}^t,\ldots,\omega_{n-1}^t). 
\end{align}
\end{proposition}

La fin de cette section est consacr\'e \`a la d\'emonstration de cette proposition. 

\begin{definition}
On note $\Omega$ l'espace des $1$-formes sur $]0,1[$ \`a valeurs dans ${\mathbb{C}}$, de la forme $\omega(z)dz=({{a_0}\over z}
+{{a_1}\over{z-1}}+f_0(z))dz$, o\`u $f_0(z)$ est une fonction lisse sur $[0,1]$ et $a_0,a_1\in{\mathbb{C}}$. On note
$\Omega_{\mathrm{reg.0}}$ et $\Omega_{\mathrm{reg.1}}$ les sous-espaces de $\Omega$ d\'efinis par les conditions respectives 
$a_0=0$ et $a_1=0$. 
\end{definition}

On a une application lin\'eaire $\Omega\to{\mathbb{C}}^2$ donn\'ee par $\omega(z)dz\mapsto(a_0,a_1)$ ; $\Omega_{\mathrm{reg.0}}$ et 
$\Omega_{\mathrm{reg.1}}$ sont alors les pr\'eimages de $0\oplus{\mathbb{C}}$ et ${\mathbb{C}}\oplus0$.  

\begin{definition}
Pour $m\geq 1$, on d\'efinit le sous-espace $(\Omega^{\otimes m})_{\mathrm{int}}$ de $\Omega^{\otimes m}$ ainsi: 
\begin{itemize}
\item si $m=1$, alors $\Omega_{\mathrm{int}}:=\Omega_{\mathrm{reg.0}}\cap\Omega_{\mathrm{reg.1}}$; 
\item si $m\geq 2$, alors $(\Omega^{\otimes m})_{\mathrm{int}}:=\Omega_{\mathrm{reg.0}}\otimes
\Omega^{\otimes m-2}\otimes\Omega_{\mathrm{reg.1}}$. 
\end{itemize}
\end{definition}
On identifie $\Omega^{\otimes m}$ \`a un sous-espace de $\Omega^m(]0,1[^m,{\mathbb{C}})$. Alors l'image de chaque \'el\'ement de 
$(\Omega^{\otimes m})_{\mathrm{int}}$ est absolument int\'egrable sur le simplexe $\Delta_m\subset]0,1[^m$. L'int\'egration sur 
$\Delta_m$ donne alors une application lin\'eaire 
$$
\int_{\Delta_m}:(\Omega^{\otimes m})_{\mathrm{int}}\to{\mathbb{C}}.
$$  
\begin{definition}
Pour $a\geq 0$, on pose $(\Omega^{\otimes a})_{{\mathrm{int.0}}}:=\Omega_{\mathrm{reg.0}}\otimes\Omega^{\otimes a-1}$
si $a>0$, et $(\Omega^{\otimes a})_{{\mathrm{int.0}}}:={\mathbb{C}}$ si $a=0$. Pour $b\geq 0$, on pose $(\Omega^{\otimes b})_{{\mathrm{int.1}}}:=\Omega^{\otimes a-1}\otimes\Omega_{\mathrm{reg.1}}$ si $b>0$, et $(\Omega^{\otimes b})_{{\mathrm{int.1}}}:={\mathbb{C}}$ si $b=0$. Enfin, pour $a,b\geq 0$, on pose $(\Omega^{\otimes a}\otimes
\Omega^{\otimes b})_{{\mathrm{int}}}:=(\Omega^{\otimes a})_{{\mathrm{int.0}}}\otimes(\Omega^{\otimes b})_{{\mathrm{int.1}}}$. 
\end{definition}

Soit $g$ un \'el\'ement de $C^\infty([0,1],{\mathbb{C}})$. Alors $dg\in\Omega$. Pour $m$ entier $\geq 1$, on d\'efinit l'application lin\'eaire
$$
{\mathrm{ins}}(dg):\bigoplus_{a+b=m-1}\Omega^{\otimes a}\otimes\Omega^{\otimes b}\to\Omega^{\otimes m}
$$
comme la somme des applications $\Omega^{\otimes a}\otimes\Omega^{\otimes b}\to\Omega^{\otimes m}$, $\alpha\otimes\beta\mapsto\alpha\otimes dg\otimes\beta$ d'insertion de $dg$ dans le produit tensoriel.  

Pour $c\geq 0$, on d\'efinit une application lin\'eaire 
$$
C^\infty([0,1],{\mathbb{C}})\otimes\Omega^{\otimes c}\to\Omega^{\otimes c}, \quad f\otimes\omega\mapsto f\cdot\omega
$$
par $f\cdot(\gamma_1\circ\cdots\circ\gamma_c):=(f\gamma_1)\circ\cdots\circ\gamma_c$ 
si $c>0$, et $f\cdot 1:=f(1)1$ pour $c=0$ ($f\gamma$ est le produit de la fonction $f$ et de la $1$-forme $\gamma$). 
On d\'efinit de m\^eme une application lin\'eaire
$$
\Omega^{\otimes c}\otimes C^\infty([0,1],{\mathbb{C}})\to\Omega^{\otimes c}, \quad \omega\otimes f\mapsto\omega\cdot f
$$
par $(\gamma_1\circ\cdots\circ\gamma_c)\cdot f:=\gamma_1\circ\cdots\circ(f\gamma_c)$ 
si $c>0$, et $1\cdot f:=f(0)1$ pour $c=0$. 

Pour $g$ dans $C^\infty([0,1],{\mathbb{C}})$, on d\'efinit alors des applications lin\'eaires 
$$
l(g),r(g):\bigoplus_{a+b=m-1}\Omega^{\otimes a}\otimes\Omega^{\otimes b}\to\Omega^{\otimes m-1}
$$
comme \'etant les sommes directes des applications $\Omega^{\otimes a}\otimes\Omega^{\otimes b}\to\Omega^{\otimes m-1}$
donn\'ees par $l(g)(\alpha\otimes\beta):=\alpha\otimes(g\cdot\beta),\quad r(g)(\alpha\otimes\beta):=(\alpha\cdot g)\otimes\beta$. 

Comme $dg\in\Omega_{\mathrm{reg.0}}\cap\Omega_{\mathrm{reg.1}}\subset\Omega$, il existe une application
$(\bigoplus_{a+b=m-1}\Omega^{\otimes a}\otimes\Omega^{\otimes b})_{{\mathrm{int}}}\to(\Omega^{\otimes m})_{{\mathrm{int}}}$, 
qui sera \'egalement not\'ee ${\mathrm{ins}}(dg)$, telle que le diagramme suivant commute 
$$
\xymatrix{
(\bigoplus_{a+b=m-1}\Omega^{\otimes a}\otimes\Omega^{\otimes b})_{{\mathrm{int}}}\ar[r]^{\ \ \ \ \ \ \ \ \ \ {\mathrm{ins}}(dg)} 
\ar@{^{(}->}[d]& 
(\Omega^{\otimes m})_{{\mathrm{int}}}\ar@{^{(}->}[d] \\ 
\bigoplus_{a+b=m-1}\Omega^{\otimes a}\otimes\Omega^{\otimes b}\ar[r]^{\ \ \ \ \ \ \ \ \ \ {\mathrm{ins}}(dg)} & \Omega^{\otimes m}}
$$
Pour tout $c\geq 0$, les applications $C^\infty([0,1],{\mathbb{C}})\otimes\Omega^{\otimes c}\to\Omega^{\otimes c}$, 
$f\otimes\omega\mapsto f\cdot\omega$ et $\Omega^{\otimes c}\otimes C^\infty([0,1],{\mathbb{C}})\to\Omega^{\otimes c}$, 
$\omega\otimes f\mapsto\omega\cdot f$, se restreignent et corestreignent en des applications $C^\infty([0,1],{\mathbb{C}})\otimes(\Omega^{\otimes c})_{{\mathrm{int.1}}}\to(\Omega^{\otimes c})_{{\mathrm{int.1}}}$ et $(\Omega^{\otimes c})_{{\mathrm{int.0}}}\otimes C^\infty([0,1],{\mathbb{C}})\to(\Omega^{\otimes c})_{{\mathrm{int.0}}}$. On en d\'eduit l'existence 
d'applications $l(g),r(g):(\bigoplus_{a+b=m-1}\Omega^{\otimes a}\otimes\Omega^{\otimes b})_{{\mathrm{int}}}\to(\Omega^{\otimes m-1})_{{\mathrm{int}}}$, telle que le diagramme
$$
\xymatrix{
(\bigoplus_{a+b=m-1}\Omega^{\otimes a}\otimes\Omega^{\otimes b})_{{\mathrm{int}}}\ar[r]^{\ \ \ \ \ \ \ \ \ \ l(g),r(g)} 
\ar@{^{(}->}[d]& 
(\Omega^{\otimes m-1})_{{\mathrm{int}}}\ar@{^{(}->}[d] \\ 
\bigoplus_{a+b=m-1}\Omega^{\otimes a}\otimes\Omega^{\otimes b}\ar[r]^{\ \ \ \ \ \ \ \ \ \ l(g),r(g)} & \Omega^{\otimes m-1}}
$$
commute. Enfin, le diagramme suivant commute 
$$
\xymatrix{
(\bigoplus_{a+b=m-1}\Omega^{\otimes a}\otimes\Omega^{\otimes b})_{{\mathrm{int}}}\ar[r]^{\ \ \ \ \ \ \ \ \ \ {\mathrm{ins}}(dg)} 
\ar[d]_{r(g)-l(g)}& 
(\Omega^{\otimes m})_{{\mathrm{int}}}\ar[d]^{\int_{\Delta_n}} \\ 
(\Omega^{\otimes m-1})_{{\mathrm{int}}}\ar[r]_{\ \ \ \ \ \ \ \ \ \ \int_{\Delta_{n-1}}} & {\mathbb{C}}}
$$
par le th\'eor\`eme de Fubini. On rassemble ces r\'esultats dans le diagramme suivant
\begin{equation}\label{diag:synth}
\xymatrix{
\bigoplus_{a+b=m-1}\Omega^{\otimes a}\otimes\Omega^{\otimes b}\ar[rr]^{{\mathrm{ins}}(dg)}\ar[dd]_{r(g)-l(g)} & & \Omega^{\otimes m} \\
& (\bigoplus_{a+b=m-1}\Omega^{\otimes a}\otimes\Omega^{\otimes b})_{{\mathrm{int}}}\ar@{_{(}->}[ul]_{{\mathrm{can}}}
\ar[r]^{\ \ \ \ \ \ \ \ \ \ {\mathrm{ins}}(dg)} 
\ar[d]_{r(g)-l(g)}& 
(\Omega^{\otimes m})_{{\mathrm{int}}}\ar[d]^{\int_{\Delta_n}}\ar@{^{(}->}[u]_{{\mathrm{can}}} \\ 
\Omega^{\otimes m-1} & (\Omega^{\otimes m-1})_{{\mathrm{int}}}\ar@{_{(}->}[l]_{{\mathrm{can}}}
\ar[r]_{\ \ \ \ \ \ \ \ \ \ \int_{\Delta_{n-1}}} & {\mathbb{C}}}
\end{equation}
o\`u chacun des carr\'es commute, et o\`u $\mathrm{can}$ sont des inclusions canoniques. 

Posons $T(\Omega):=\oplus_{a\geq 0}\Omega^{\otimes a}$. Les sommes sur $m\geq 0$ des applications lin\'eaires
$$
l(g),r(g):\bigoplus_{a,b|a+b=m}\Omega^{\otimes a}\otimes\Omega^{\otimes b}\to\Omega^{\otimes m}
$$ 
sont des applications lin\'eaires $l(g),r(g):T(\Omega)^{\otimes2}\to T(\Omega)$.  

Muni du produit de battage $\sh$ et du coproduit de 
d\'econcat\'enation $\Delta_\sh$ donn\'e par $\Delta_\sh(\alpha_1\circ\cdots\circ\alpha_n):=\sum_{k=0}^n (\alpha_1\circ\cdots\circ\alpha_k)\otimes(\alpha_{k+1}\circ\cdots
\circ\alpha_n)$, l'espace $T(\Omega)$ est une big\`ebre. Le produit $\sh$ \'etant commutatif, les repr\'esentations r\'eguli\`eres \`a gauche et \`a droite de $T(\Omega)$ sont isomorphes ; ceci donne un $T(\Omega)$-module $T(\Omega)$. En utilisant le coproduit, on munit l'espace $T(\Omega)^{\otimes 2}$ d'une structure de $T(\Omega)$-module. On a alors : 

\begin{lemma}
L'application lin\'eaire $r(g)-l(g) :T(\Omega)^{\otimes 2}\to T(\Omega)$ est un morphisme de $T(\Omega)$-modules. 
En d'autres termes, on a 
\begin{equation}\label{id:shuffle}
(r(g)-l(g))\big((\alpha^{(1)}\sh\beta)\otimes(\alpha^{(2)}\sh\gamma)\big)
=\alpha\sh\big((r(g)-l(g))(\beta\otimes\gamma)\big)
\end{equation}
(\'egalit\'e dans $T(\Omega)$) pour tous $\alpha,\beta,\gamma$ dans $T(\Omega)$, dans laquelle on note 
$\Delta_\sh(\alpha)=\alpha^{(1)}\otimes\alpha^{(2)}$. 
\end{lemma} 

{\em D\'emonstration.} Notons $a(g):=r(g)-l(g)$. Si $\beta=\gamma=1$, le membre de gauche est \'egal \`a $a(g)(\alpha^{(1)}\otimes\alpha^{(2)})$ qui est \'egal, apr\`es simplifications, \`a $\alpha\circ(g\cdot 1)-(1\cdot g)\circ\alpha=(g(1)-g(0))\alpha$. D'autre part, le membre de droite est \'egal \`a $\alpha\sh\big(a(g)(1\otimes 1)\big)=\alpha\sh\big((g\cdot 1)-(1\cdot g)\big)=(g(1)-g(0))\alpha$. {\it L'identit\'e (\ref{id:shuffle}) est donc satisfaite si $\beta=\gamma=1$.} 

Si $\beta\in\Omega$ et $\gamma=1$, alors le membre de gauche est \'egal \`a  $a(g)((\alpha^{(1)}\sh\beta)\otimes\alpha^{(2)})$, ce 
qui d'apr\`es l'identit\'e 
\begin{equation}\label{égalité :utile}
\forall \alpha\in T(\Omega), \forall \beta\in\Omega, \quad \alpha\sh\beta=\alpha^{(1)}\circ\beta\circ\alpha^{(2)},
\end{equation} 
est \'egal \`a $a(g)((\alpha^{(1)}\circ\beta\circ\alpha^{(2)})\otimes\alpha^{(3)})$. Apr\`es simplification, ce dernier terme est \'egal \`a 
$g(1) \alpha^{(1)}\circ \beta\circ\alpha^{(2)}-\alpha^{(1)}\circ(g\beta)\circ\alpha^{(2)}$.  D'autre part, le membre de droite est \'egal \`a $\alpha\sh\big(a(g)(\beta\otimes 1)\big)$. On a $a(g)(\beta\otimes 1)=g(1)\beta-g\beta$, donc compte tenu de (\ref{égalité :utile}), 
le membre de droite est \'egal \`a $\alpha^{(1)}\circ(g(1)\beta-g\beta)\circ\alpha^{(2)}$ {\it L'identit\'e (\ref{id:shuffle}) est donc satisfaite si $\beta\in\Omega$ et $\gamma=1$.} On montre de m\^eme que {\it l'identit\'e (\ref{id:shuffle}) est satisfaite si $\beta=1$ et $\gamma\in\Omega$.}

Si $\beta$ et $\gamma$ appartiennent \`a $\Omega$, le membre de gauche est \'egal \`a 
$a(g)\big((\alpha^{(1)}\sh\beta)\otimes(\alpha^{(2)}\sh\gamma)\big)$, ce qui 
d'apr\`es (\ref{égalité :utile}) est \'egal \`a $a(g)\big((\alpha^{(1)}\circ\beta\circ
\alpha^{(2)}) \otimes(\alpha^{(3)}\circ\gamma\circ\alpha^{(4)})\big)$, terme qui compte 
tenu de l'\'egalit\'e suivante (dans laquelle $T(\Omega)_+=\oplus_{a>0}\Omega^{\otimes a}$)
$$
\forall\alpha,\beta\in T(\Omega),\forall\alpha’,\beta’\in T(\Omega)_+, \quad a(g)
((\alpha\circ\alpha')\otimes(\beta'\circ\beta))=\alpha\circ a(g)(\alpha'\otimes\beta')\circ\beta
$$
est \'egal \`a $\alpha^{(1)}\circ a(g)\big( (\beta\circ\alpha^{(2)}) \otimes(\alpha^{(3)}
\circ\gamma)\big) \circ\alpha^{(4)}$. 
L'\'egalit\'e suivante
$$
\forall\beta,\gamma\in\Omega,\forall\alpha\in T(\Omega), \quad
a(g)\big((\beta\circ\alpha^{(1)})\otimes(\alpha^{(2)}\circ\gamma)\big) =
\beta\circ\alpha\circ(g\gamma)-(g\beta)\circ\alpha\circ\gamma
$$ 
implique alors que ce dernier terme est \'egal \`a $\alpha^{(1)}\circ \big(
\beta\circ\alpha^{(2)}\circ(g\gamma) - (g\beta)\circ\alpha^{(2)}\circ\gamma \big) \circ\alpha^{(3)}$. 
On a l'identit\'e suivante 
\begin{equation}\label{star}
\forall\beta,\gamma\in\Omega,\forall\alpha\in T(\Omega), \quad \alpha\sh
(\beta\circ\gamma)=\alpha^{(1)}\circ\beta\circ\alpha^{(2)}\circ\gamma\circ\alpha^{(3)},  
\end{equation}
qui permet d'exprimer le dernier terme comme $\alpha\sh\big(\beta\circ(g\gamma)-(g\beta)\circ\gamma\big)$, 
qui est donc $\alpha\sh a(g)(\beta\otimes\gamma)$, et donc \'egal au membre de droite. 
{\it L'identit\'e (\ref{id:shuffle}) est donc satisfaite si $\beta=1$ et $\gamma\in\Omega$.}

Les deux membres de l'identit\'e (\ref{id:shuffle}) sont les valeurs en $\alpha\otimes\beta\otimes\gamma$ de 
deux applications lin\'eaires $gch,dt :T(\Omega)\otimes T(\Omega)^{\otimes 2}\to T(\Omega)$. 
La d\'ecomposition 
$$
T(\Omega)^{\otimes 2}=T(\Omega)_+^{\otimes 2}\oplus (T(\Omega)_+\otimes{\mathbb{C}})
\oplus({\mathbb{C}}\otimes T(\Omega)_+)\oplus{\mathbb{C}}^{\otimes 2}
$$ 
induit par tensorisation avec $T(\Omega)$ une d\'ecomposition 
de la source des applications $gch,dt$. On a par ailleurs montr\'e l'\'egalit\'e des restrictions des applications 
$gch$ et $dt$ aux produits tensoriels de $T(\Omega)$ avec les sous-espaces $\Omega^{\otimes 2}\subset 
T(\Omega)_+^{\otimes 2}$, $(\Omega\otimes{\mathbb{C}})\subset (T(\Omega)\otimes{\mathbb{C}})$, 
$({\mathbb{C}}\otimes T(\Omega))\subset ({\mathbb{C}}\otimes T(\Omega))$, et ${\mathbb{C}}^{\otimes 2}\subset
{\mathbb{C}}^{\otimes 2}$. L'\'egalit\'e de ces applications lin\'eaires sur les trois premiers espaces est 
alors une cons\'equence des identit\'es
$$
\forall\tilde\beta\in T(\Omega),\forall\beta_0\in T(\Omega)_+,\quad a(g)\big((\tilde\beta\circ\beta_0)\otimes 1\big)
=\tilde\beta\circ a(g)(\beta_0\otimes 1), 
$$
$$
\forall\tilde\gamma\in T(\Omega),\forall\gamma_0\in T(\Omega)_+,\quad a(g)\big(1\otimes 
(\gamma_0\circ\tilde\gamma)\big)=a(g)(1\otimes\gamma_0)\circ\tilde\gamma, 
$$
\begin{equation}\label{help}
\forall\tilde\beta,\tilde\gamma\in T(\Omega),\forall\beta,\gamma\in T(\Omega)_+,\quad a(g)
\big((\tilde\beta\circ\beta_0)\otimes(\gamma_0\circ\tilde\gamma)\big)=\tilde\beta\circ 
a(g)(\beta_0\otimes\gamma_0)\circ\tilde\gamma,  
\end{equation}
et de l'identit\'e 
\begin{equation}\label{shuffle:conc}
\alpha\sh(\beta\circ\gamma)=(\alpha^{(1)}\sh\beta)\circ(\alpha^{(2)}\sh\gamma)
\end{equation} 
pour $\alpha,\beta,\gamma\in T(\Omega)$. 

Montrons par exemple l'identit\'e (\ref{id:shuffle}) dans le cas du produit tensoriel de $T(\Omega)$ avec 
$T(\Omega)_+^{\otimes 2}$. Par lin\'earit\'e, on suppose $\beta\otimes\gamma\in T(\Omega)_+^{\otimes 2}$
de la forme $(\tilde\beta\circ\beta_0)\otimes(\gamma_0\circ\tilde\gamma)$ o\`u $\tilde\beta,\tilde\gamma\in T(\Omega)$ 
et $\beta_0,\gamma_0\in\Omega$. Alors, si $\alpha\in T(\Omega)$, 
\begin{align*}
& a(g)\big((\alpha^{(1)}\sh\beta)\otimes(\alpha^{(2)}\sh\gamma)\big)
=a(g)\Big(\big(\alpha^{(1)}\sh(\tilde\beta\circ\beta_0)\big)\otimes\big(\alpha^{(2)}\sh(\gamma_0\circ\tilde\gamma)\big)\Big) 
\\ & =a(g)\Big( \big((\alpha^{(1)}\sh\tilde\beta)\circ(\alpha^{(2)}\sh\beta_0)\big)
\otimes\big((\alpha^{(3)}\sh\gamma_0)\circ(\alpha^{(4)}\sh\tilde\gamma)\big)\Big) \quad
\text{ (en utilisant (\ref{star}))}
\\ & =(\alpha^{(1)}\sh\tilde\beta)\circ a(g)\big((\alpha^{(2)}\sh\beta_0)
\otimes(\alpha^{(3)}\sh\gamma_0)\big)\circ(\alpha^{(4)}\sh\tilde\gamma)\quad
\text{ (en utilisant (\ref{help}))}
\\ & =(\alpha^{(1)}\sh\tilde\beta)\circ 
\Big( \alpha^{(2)}\sh a(g)\big(\beta_0\otimes\gamma_0\big)\Big)
\circ(\alpha^{(3)}\sh\tilde\gamma)\quad
\text{ (en utilisant (\ref{id:shuffle}) pour $\beta_0,\gamma_0\in\Omega$)}
\\ & =\alpha\sh\Big(\tilde\beta\circ a(g)\big(\beta_0\otimes\gamma_0\big)
\circ\tilde\gamma\Big)\quad \text{ (en utilisant (\ref{shuffle:conc}))}
\\ & =\alpha\sh a(g)\big((\tilde\beta\circ\beta_0)\otimes(\gamma_0\circ\tilde\gamma)\big)
\quad\text{ (en utilisant (\ref{help}))}
\\ & =\alpha\sh  a(g)(\beta\otimes\gamma). 
\end{align*}
\hfill\qed\medskip
 
Compte tenu de ce que $(T(\Omega),\sh,\Delta_\sh)$ est une big\`ebre, (\ref{id:shuffle}) implique : 
\begin{equation}\label{id:shuffle:bis}
\forall\alpha,\beta,\gamma,\delta\in T(\Omega), \quad 
(r(g)-l(g))\big((\alpha^{(1)}\sh\beta\sh\delta^{(1)})\otimes(\alpha^{(2)}\sh\gamma\sh\delta^{(2)})\big)
=\alpha\sh\big((r(g)-l(g))(\beta\otimes\gamma)\big)\sh\delta
\end{equation}
(\'egalit\'e dans $T(\Omega)$). 

Pour $\breve\omega_1,\ldots,\breve\omega_n\in\Omega_{1,1}$, on pose 
$$
{\mathrm{int}}(\breve\omega_1,\ldots,\breve\omega_n):=\sum_{a,b|a+b\leq n}(-1)^{a+b}(d{\mathrm{ln}}(z))^{\circ a}
\sh(\breve\omega_{a+1}\circ\cdots\circ\breve\omega_{n-b})\sh (d{\mathrm{ln}}(1-z))^{\circ b}\in \Omega^{\otimes n}. 
$$
On a vu en section \ref{subsect:IIreg:holon} que 
${\mathrm{int}}(\breve\omega_1,\ldots,\breve\omega_n)=\sum_{{\epsilon,\eta\in\{0,1\}}
\atop{k,l|k+l\leq n-2}}{\mathrm{int}}_{k,l}^{\epsilon,\eta}(\breve\omega_1,\ldots,\breve\omega_n)$, o\`u
\begin{align*}
&{\mathrm{int}}_{k,l}^{0,0}(\breve\omega_1,\ldots,\breve\omega_n)\\ & =
 (\breve\omega_{k+1}-d{\mathrm{log}}(z))\circ\Big((-d{\mathrm{log}}(z))^{\circ k}\sh(-d{\mathrm{log}}(1-z))^{\circ l}\sh
(\breve\omega_{k+2}\circ\cdots\circ\breve\omega_{n-l-1})\Big)\circ(\breve\omega_{n-l}-d{\mathrm{log}}(1-z)), 
\end{align*}
\begin{align*}
&{\mathrm{int}}_{k,l}^{1,0}(\breve\omega_1,\ldots,\breve\omega_n)\\ & =
 (\breve\omega_{k+2}-d{\mathrm{log}}(z))\circ\Big((-d{\mathrm{log}}(z))^{\circ k}\sh(-d{\mathrm{log}}(1-z))^{\circ l}\sh
(\breve\omega_{k+3}\circ\cdots\circ\breve\omega_{n-l})\Big)\circ(-d{\mathrm{log}}(z)), 
\end{align*}
\begin{align*}
&{\mathrm{int}}_{k,l}^{0,1}(\breve\omega_1,\ldots,\breve\omega_n)\\ & =
 (-d{\mathrm{log}}(1-z))\circ\Big((-d{\mathrm{log}}(z))^{\circ k}\sh(-d{\mathrm{log}}(1-z))^{\circ l}\sh
(\breve\omega_{k+1}\circ\cdots\circ\breve\omega_{n-l-2})\Big)\circ(\breve\omega_{n-l-1}-d{\mathrm{log}}(1-z)), 
\end{align*}
\begin{align*}
&{\mathrm{int}}_{k,l}^{0,0}(\breve\omega_1,\ldots,\breve\omega_n)\\ & =
 (-d{\mathrm{log}}(1-z))\circ\Big((-d{\mathrm{log}}(z))^{\circ k}\sh(-d{\mathrm{log}}(1-z))^{\circ l}\sh
(\breve\omega_{k+2}\circ\cdots\circ\breve\omega_{n-l-1})\Big)\circ(-d{\mathrm{log}}(z)).  
\end{align*}
Chaque ${\mathrm{int}}_{k,l}^{\epsilon,\eta}(\breve\omega_1,\ldots,\breve\omega_n)$ appartient \`a 
$(\Omega^{\otimes n})_{{\mathrm{int}}}$, donc ${\mathrm{int}}(\breve\omega_1,\ldots,\breve\omega_n)$ aussi. On a aussi pos\'e
$$
I_{[0,1]}^{reg}(\breve\omega_1,\ldots,\breve\omega_n)=\int_{\Delta_n}{\mathrm{int}}(\breve\omega_1,\ldots,\breve\omega_n). 
$$
Soit $\breve g_1,\ldots,\breve g_n\in C^\infty([0,1],{\mathbb{C}})$ tels que $\breve g_i(0)=\breve g_i(1)$ ; pour $n\in\{1,\ldots,n\}$, 
on d\'efinit l'\'el\'ement suivant $\Omega^{\otimes n}$
\begin{align*}
& \delta_c{\mathrm{int}}(\breve\omega_1,\ldots,\breve\omega_n):=
\sum_{{a,b|a+b\leq n}\atop{a+1\leq c\leq n-b}}(-1)^{a+b}(-1)^{a+b}(d{\mathrm{ln}}(z))^{\circ a}\sh
(\breve\omega_{a+1}\circ\cdots\circ d\breve g_c\circ\cdots\circ\breve\omega_{n-b})\sh(d{\mathrm{ln}}(1-z))^{\circ b}. 
\end{align*}
On d\'efinit $\delta_c{\mathrm{int}}_{k,l}^{\epsilon,\eta}(\breve\omega_1,\ldots,\breve\omega_n)$ comme \'etant 0 si 
$c\notin\{k+1,\ldots,n-l\}$, et comme \'etant le r\'esultat du remplacement du terme comprenant $\breve\omega_c$
(qui peut \^etre $\breve\omega_c$, $\breve\omega_c-d{\mathrm{log}}(z)$ ou $\breve\omega_c-d{\mathrm{log}}(1-z)$) par $d\breve g_c$. 
Par exemple, on a 
\begin{align*}
& \delta_c{\mathrm{int}}_{k,l}^{0,0}(\breve\omega_1,\ldots,\breve\omega_n):=
\\ & \delta_{k+1,c}d\breve g_c\circ\Big((-d{\mathrm{log}}(z))^{\circ k}\sh(-d{\mathrm{log}}(1-z))^{\circ l}\sh
(\breve\omega_{k+2}\circ\cdots\circ\breve\omega_{n-l-1})\Big)\circ(\breve\omega_{n-l}-d{\mathrm{log}}(1-z))
\\ & +\delta_{k+2\leq c\leq n-l-1}
(\breve\omega_{k+1}-d{\mathrm{log}}(z))\circ\\ & \circ\Big((-d{\mathrm{log}}(z))^{\circ k}\sh(-d{\mathrm{log}}(1-z))^{\circ l}\sh
(\breve\omega_{k+2}\circ\cdots\circ d\breve g_c\circ\cdots\circ\breve\omega_{n-l-1})\Big)\circ(\breve\omega_{n-l}-d{\mathrm{log}}(1-z))
\\ & +\delta_{n-l,c}
(\breve\omega_{k+1}-d{\mathrm{log}}(z))\circ\Big((-d{\mathrm{log}}(z))^{\circ k}\sh(-d{\mathrm{log}}(1-z))^{\circ l}\sh
(\breve\omega_{k+2}\circ\cdots\circ\breve\omega_{n-l-1})\Big)\circ d\breve g_c 
\end{align*}
(o\`u $\delta_{k+2\leq c\leq n-l-1}=1$ si $k+2\leq c\leq n-l-1$, et $=0$ sinon). 

Comme $d\breve g_c\in\Omega_{{\mathrm{reg.0}}}\cap\Omega_{{\mathrm{reg.1}}}$, on a $\delta_c{\mathrm{int}}_{k,l}^{0,0}(\breve\omega_1,
\ldots,\breve\omega_n)\in(\Omega^{\otimes n})_{\mathrm{int}}$ pour tout $(k,l)$. De m\^eme, on a pour 
$(\epsilon,\eta)\in\{(1,0),(0,1),(1,1)\}$ et tout $(k,l)$, $\delta_c{\mathrm{int}}_{k,l}^{\epsilon,\eta}(\breve\omega_1,\ldots,\breve\omega_n)\in(\Omega^{\otimes n})_{\mathrm{int}}$. On en d\'eduit
$$
\delta_c{\mathrm{int}}(\breve\omega_1,\ldots,\breve\omega_n)\in(\Omega^{\otimes n})_{\mathrm{int}}. 
$$ 
On d\'efinit l'\'el\'ement $\mathrm{rel}_c\in\bigoplus_{a,b|a+b=n-1}\Omega^{\otimes a}\otimes\Omega^{\otimes b}$ 
par 
\begin{align*}
\mathrm{rel}_c:=\sum_{{a,b|a+b\leq n}\atop{a+1\leq c\leq n-b}}\sum_{{a=a'+a''}\atop{b=b'+b''}}
& (-1)^{a+b}\Big((d\mathrm{ln}(z))^{\circ a'}\sh(\breve\omega_{a+1}\circ\cdots\circ\breve\omega_{c-1})\sh(d\mathrm{ln}(1-z))^{
\circ b'}\Big)\otimes\\ &\otimes\Big((d\mathrm{ln}(z))^{\circ a''}\sh(\breve\omega_{c+1}\circ\cdots\circ\breve\omega_{n-b})
\sh(d\mathrm{ln}(1-z))^{\circ b''}\Big)
\end{align*}
On d\'efinit $(\mathrm{rel}_c)_{k,l}^{0,0}\in\bigoplus_{a,b|a+b=n-1}(\Omega^{\otimes a}\otimes\Omega^{\otimes b})_{\mathrm{int}}$
par 
\begin{align*}
&(\mathrm{rel}_c)_{k,l}^{0,0}:=
\delta_{k+1,c}\cdot 1\otimes\Big(\big((-d{\mathrm{log}}(z))^{\circ k}\sh(-d{\mathrm{log}}(1-z))^{\circ l}\sh
(\breve\omega_{k+2}\circ\cdots\circ\breve\omega_{n-l-1})\big)\circ(\breve\omega_{n-l}-d{\mathrm{log}}(1-z))\Big)
\\ & +\delta_{k+2\leq c\leq n-l-1}
(\breve\omega_{k+1}-d{\mathrm{log}}(z))\circ\\ & \circ\Big((-d{\mathrm{log}}(z))^{\circ k}\sh(-d{\mathrm{log}}(1-z))^{\circ l}\sh
(\breve\omega_{k+2}\circ\cdots\circ d\breve g_c\circ\cdots\circ\breve\omega_{n-l-1})\Big)\circ(\breve\omega_{n-l}-d{\mathrm{log}}(1-z))
\\ & +\sum_{{k',k''|k=k'+k''}\atop{l',l''|l=l'+l''}}\delta_{k+2\leq c\leq n-l-1}\Big((\breve\omega_{k+1}-d{\mathrm{log}}(z))\circ
\big((-d\mathrm{ln}(z))^{\circ k'}\sh(-d\mathrm{ln}(1-z))^{\circ l'}\sh(\breve\omega_{k+2}\circ\cdots\circ\breve\omega_{c-1})\big)\Big)\otimes
\\ & \otimes
\Big(\big((\breve\omega_{c+1}\circ\cdots\circ\breve\omega_{n-l-1})\sh(-d\mathrm{ln}(z))^{\circ k''}\sh(-d\mathrm{ln}(1-z))^{\circ l''}\big)\circ(\breve\omega_{n-l}-d{\mathrm{log}}(1-z))\Big)
\\ & +\delta_{n-l,c}\Big(
(\breve\omega_{k+1}-d{\mathrm{log}}(z))\circ\big((-d{\mathrm{log}}(z))^{\circ k}\sh(-d{\mathrm{log}}(1-z))^{\circ l}\sh
(\breve\omega_{k+2}\circ\cdots\circ\breve\omega_{n-l-1})\big)\Big)\otimes 1
\end{align*}
et on d\'efinit de fa\c con analogue les autres $(\mathrm{rel}_c)_{k,l}^{\epsilon,\eta}$ dans 
$\bigoplus_{a,b|a+b=n-1}(\Omega^{\otimes a}\otimes\Omega^{\otimes b})_{\mathrm{int}}$. Alors $\mathrm{rel}_c=\sum_{\epsilon,\eta≤\in\{0,1\}}\sum_{k,l|k+l\leq n-2}(\mathrm{rel}_c)_{k,l}^{\epsilon,\eta}$, donc 
$$
\mathrm{rel}_c\in\bigoplus_{a,b|a+b=n-1}(\Omega^{\otimes a}\otimes\Omega^{\otimes b})_{\mathrm{int}}. 
$$
Compte tenu de (\ref{id:shuffle:bis}), on a 
\begin{align*}
(r(\breve g_c)-l(\breve g_c))(\mathrm{rel}_c)=\sum_{{a,b|a+b\leq n}\atop{a+1\leq c\leq n-b}}&(-1)^{a+b}
(d{\mathrm{log}}(z))^{\circ a}\sh
(d{\mathrm{log}}(1-z))^{\circ b}\sh\\ & \sh(r(\breve g_c)-l(\breve g_c))((\breve\omega_{a+1}\circ\cdots\circ\breve\omega_{c-1})\otimes
(\breve\omega_{c+1}\circ\cdots\circ\breve\omega_{n-b}))
\end{align*}
On a donc 
\begin{align*}
& \sum_{c=1}^n (l(\breve g_c)-r(\breve g_c))(\mathrm{rel}_c)=\sum_{a,b|a+b\leq n}(-1)^{a+b}(d{\mathrm{log}}(z))^{\circ a}
\sh(d{\mathrm{log}}(1-z))^{\circ b}\sh
\\ & \sh\Big(-\breve g_{a+1}(0)(\breve\omega_{a+2}\circ\cdots\circ\breve\omega_{n-b})+(\breve\omega_{a+1}\circ\cdots\circ
\breve\omega_{n-b-1})g_{n-b}(0)
\\ & +\sum_{i=a+1}^{n-b-1}\breve\omega_{a+1}\circ\cdots\circ(\breve g_i\breve\omega_{i+1}-\breve g_{i+1}\breve\omega_i)\circ\cdots\circ\breve\omega_{n-b}\Big). 
\end{align*}
Supposons que pour $i=1,\ldots,n-1$, on dispose de $\breve\psi_{i,i+1}\in\Omega_{1,1}$ tel que 
$\breve g_i\breve\omega_{i+1}-\breve g_{i+1}\breve\omega_i=(\breve g_i(0)-\breve g_{i+1}(0))\breve\psi_{i,i+1}$. Alors 
\begin{align}\label{égalité:interm:1}
& \nonumber \sum_{c=1}^n (l(\breve g_c)-r(\breve g_c))(\mathrm{rel}_c)=\sum_{a,b|a+b\leq n}(-1)^{a+b}(d{\mathrm{log}}(z))^{\circ a}
\sh(d{\mathrm{log}}(1-z))^{\circ b}\sh
\\ & \nonumber \sh\Big(-\breve g_{a+1}(0)(\breve\omega_{a+2}\circ\cdots\circ\breve\omega_{n-b})+(\breve\omega_{a+1}
\circ\cdots\circ\breve\omega_{n-b-1})\breve g_{n-b}(0)
\\ & +\sum_{i=a+1}^{n-b-1}(\breve g_i(0)-\breve g_{i+1}(0))\breve\omega_{a+1}\circ\cdots\circ\breve\psi_{i,i+1}\circ\cdots\circ
\breve\omega_{n-b}\Big) 
\end{align}
(\'egalit\'e dans $\Omega^{\otimes n-1}$)
D'autre part, en d\'ecomposant le d\'eveloppement de $\mathrm{int}(\breve\omega_1,\ldots,\breve\psi_{i,i+1},\ldots,\breve\omega_n)$ 
selon les valeurs de $a,b$, on a le d\'eveloppement suivant dans $\Omega^{\otimes n-1}$ 
\begin{align}\label{égalité:interm:2}
&\nonumber -\breve g_1(0)\mathrm{int}(\breve\omega_2,\ldots,\breve\omega_n)+\breve g_n(0)\mathrm{int}(\breve\omega_1,\ldots,
\breve\omega_{n-1})
+\sum_{i=1}^{n-1}(\breve g_i(0)-\breve g_{i+1}(0))\mathrm{int}(\breve\omega_1,\ldots,\breve\psi_{i,i+1},\ldots,\breve\omega_n)
\\ & \nonumber
=-\breve g_1(0)\sum_{a,b|a+b\leq n-1}(-d\mathrm{ln}(z))^{\circ a}\sh(-d\mathrm{ln}(1-z))^{\circ b}\sh
(\breve\omega_{a+2}\circ\cdots\circ\breve\omega_{n-b})
\\ & \nonumber+\breve g_n(0)\sum_{a,b|a+b\leq n-1}(-d\mathrm{ln}(z))^{\circ a}\sh(-d\mathrm{ln}(1-z))^{\circ b}\sh
(\breve\omega_{a+1}\circ\cdots\circ\breve\omega_{n-b-1})
\\ & \nonumber+\sum_{i=1}^{n-1}(\breve g_i(0)-\breve g_{i+1}(0))\sum_{{a,b|a+b\leq n-1}\atop{a\leq i-1,b\leq n-i-1}}
(-d\mathrm{ln}(z))^{\circ a}\sh(-d\mathrm{ln}(1-z))^{\circ b}\sh
(\breve\omega_{a+1}\circ\cdots\circ\breve\psi_{i,i+1}\circ\cdots\circ\breve\omega_{n-b})
\\ & \nonumber+\sum_{i=1}^{n-1}(\breve g_i(0)-\breve g_{i+1}(0))\sum_{{a,b|a+b\leq n-1}\atop{a\geq i-1}}
(-d\mathrm{ln}(z))^{\circ a}\sh(-d\mathrm{ln}(1-z))^{\circ b}\sh
(\breve\omega_{a+2}\circ\cdots\circ\breve\omega_{n-b})
\\ & +\sum_{i=1}^{n-1}(\breve g_i(0)-\breve g_{i+1}(0))\sum_{{a,b|a+b\leq n-1}\atop{b\geq n-i}}
(-d\mathrm{ln}(z))^{\circ a}\sh(-d\mathrm{ln}(1-z))^{\circ b}\sh
(\breve\omega_{a+1}\circ\cdots\circ\breve\omega_{n-b-1})
\end{align}
Apr\`es inversion des signes somme dans le quatri\`eme terme du membre de droite de (\ref{égalité:interm:2}), on 
voit que la somme de ce terme et du premier terme de ce membre de droite est \'egale au premier terme du membre de 
droite de (\ref{égalité:interm:1}). De m\^eme, la somme des deuxi\`eme et cinqui\`eme termes du membre de droite 
de (\ref{égalité:interm:2}) est \'egale au deuxi\`eme terme du membre de droite de (\ref{égalité:interm:1}). Enfin, les trois\`emes
termes des membres de droite de (\ref{égalité:interm:1}) et (\ref{égalité:interm:2}) sont \'egaux. On a donc l'\'egalit\'e 
suivante 
\begin{align}\label{last:eq}
& \sum_{c=1}^n (l(\breve g_c)-r(\breve g_c))(\mathrm{rel}_c)
\\ & \nonumber = -\breve g_1(0)\mathrm{int}(\breve\omega_2,\ldots,\breve\omega_n)+\breve g_n(0)\mathrm{int}(\breve\omega_1,
\ldots,\breve\omega_{n-1})
+\sum_{i=1}^{n-1}(\breve g_i(0)-\breve g_{i+1}(0))\mathrm{int}(\breve\omega_1,\ldots,\breve\psi_{i,i+1},\ldots,\breve\omega_n)
\end{align}
dans $\Omega^{\otimes n-1}$, et donc dans $(\Omega^{\otimes n-1})_{\mathrm{int}}$. 

Rappelons que pour $c\in\{1,\ldots,n\}$, $\mathrm{rel}_c\in(\bigoplus_{a+b=n-1}\Omega^{\otimes a}\otimes\Omega^{\otimes b})_{\mathrm{int}}$ est tel que $\mathrm{ins}(d\breve g_c)(\mathrm{rel}_c)=\delta_c\mathrm{int}(\breve\omega_1,\ldots,\breve\omega_n)$. 
On en d\'eduit (voir (\ref{diag:synth})) que $(r(\breve g_c)-l(\breve g_c))(\mathrm{rel}_c)\in(\Omega^{\otimes n-1})_{\mathrm{int}}$ et 
$\int_{\Delta_n}\delta_c\mathrm{int}(\breve\omega_1,\ldots,\breve\omega_n)=\int_{\Delta_{n-1}}(r(\breve g_c)-l(\breve g_c))(\mathrm{rel}_c)$. 
En sommant sur $c=1,\ldots,n$, et en utilisant (\ref{last:eq}), on en d\'eduit 
\begin{equation}\label{truc}
\int_{\Delta_n}\delta_c\mathrm{int}(\breve\omega_1,\ldots,\breve\omega_n)=\int_{\Delta_{n-1}}(\text{membre de droite de (\ref{last:eq})}).
\end{equation}

On se place maintenant dans le cadre de l'\'enonc\'e de la proposition \ref{prop:3:1}. Fixons $t\in I$ et posons $\breve\omega_i:=\omega_i^t$, 
$\breve g_i:=g_i^t$, $\breve\psi_{i,i+1}:=\psi_{i,i+1}^t$. Alors les hypoth\`eses sur les $\breve\omega_i$, $\breve g_i$, $\breve\psi_{i,i+1}$ sont satisfaites, et $\delta_c\mathrm{int}(\breve\omega_1,\ldots,\breve\omega_n)=(d/dt)\mathrm{int}(\omega_1^t,\ldots,\omega_n^t)$, donc le membre de gauche de (\ref{truc}) s'identifie au membre de gauche de (\ref{eq:diff}). De m\^eme, le membre de droite de (\ref{truc}) s'identifie au membre de droite de (\ref{eq:diff}). Ceci montre (\ref{eq:diff}) et termine la d\'emonstration de la proposition \ref{prop:3:1}. 

\subsection{Syst\`emes diff\'erentiels pour les analogues elliptiques des nombres multiz\'etas}\label{sec:3:2}

Soient $x_1,\ldots,x_n$ des nombres complexes (ou des variables formelles proches de 0). 

\begin{definition}\label{def:ell}
Pour $i\in\{1,\ldots,n\}$ et $(\tau,z)\in{\mathfrak{H}}\times]0,1[$, on pose 
$$
\omega_i(\tau,z)dz :=\sigma_{x_i}^\tau(z)dz, \quad 
g_i(\tau,z)d\tau :={\partial\over{\partial x_i}}({1\over{2\pi\mathrm{i}}}\sigma_{x_i}^{\tau}(z))d\tau, 
$$ 
et pour $i\in\{1,\ldots,n-1\}$ et $(\tau,z)\in{\mathfrak{H}}\times]0,1[$, 
$$
\psi_{i,i+1}(\tau,z)dz :=\sigma_{x_i+x_{i+1}}^\tau(z)dz. 
$$
\end{definition}

On a vu que pour chaque couple $(\tau,x)$, $\sigma_x^\tau(z)dz$ est dans 
$\Omega_{1,1}$, ce qui implique que les $\omega_i$ ($i=1,\ldots,n$) et les $\psi_{i,i+1}$ 
($i=1,\ldots,n-1$) sont dans $\mathcal{F}$. De plus, on a le d\'eveloppement 
\begin{align*}
& {1\over{2\pi\mathrm{i}}}{1\over{\theta_\tau(z)}}{{\theta_\tau(z+x)}\over{\theta_\tau(x)}}
={1\over{2\pi\mathrm{i}}}({1\over z}+O(z)){{\theta_\tau(z+x)}\over{\theta_\tau(x)}}  
\\ &={1\over{2\pi\mathrm{i}\theta_\tau(x)}}({{\theta_\tau(x)}\over{z}}+\theta'_\tau(x)+O(z)) 
={1\over{2\pi\mathrm{i}}}({1\over z}+{{\theta'_\tau}\over{\theta_\tau}}(x)+O(z))
\end{align*} 
d'o\`u le fait que \`a $(\tau,x)$ fix\'e, $z\mapsto g_x(\tau,z)$ est lisse en $0$, avec 
$$
g_x(\tau,z)={\partial\over{\partial x}}
({1\over{2\pi\mathrm{i}}}{{\theta’_\tau}\over{\theta_\tau}}(x))+O(z), 
$$ 
d'o\`u l'on d\'eduit 
\begin{equation}\label{formule:g:i}
g_x(\tau,0)={\partial\over{\partial x}}({1\over{2\pi\mathrm{i}}}{{\theta’_\tau}\over{\theta_\tau}}(x)).
\end{equation}
L'invariance de $z\mapsto\theta_\tau(z)$ sous $z\mapsto 1-z$ implique que $g_z(\tau,z)$ est 
\'egalement lisse en $1$, avec $g_x(\tau,0)=g_x(\tau,1)$. 

Dans \cite{CEE}, trois lignes apr\`es l'\'equation (14), on montre l'\'egalit\'e 
$$
\partial_{\tau}({{\theta_\tau(z+x)}\over{\theta_\tau(z)\theta_\tau(x)}}) = 
{1\over{2\pi\on{i}}}\partial_{z}\partial_{x}({{\theta_\tau(z+x)}\over{\theta_\tau(z)\theta_\tau(x)}})
$$
qui implique imm\'ediatement 
$$
{{\partial\omega_i}\over{\partial\tau}}(\tau,z)={{\partial g_i}\over{\partial z}}(\tau,z)
$$
pour $i\in\{1,\ldots,n\}$. 

La fonction de Weierstrass est d\'efinie par $\wp_\tau(z) = \sum'_{a\in \ZZ+\tau\ZZ}
((z+a)^{-2}- a^{-2})$, o\`u $\sum'$ signifie que le terme
$a^{-2}$ n'est pas pris en compte lorsque $a=0$. 
On pose alors $$\tilde\wp_\tau(z):= \wp_\tau(z)+G_{2}(\tau) = \sum_{m\in\ZZ}
(\sum_{n}{}^{'}(z+n+m\tau)^{-2}).$$ 

\begin{lemma} \label{lemme:weier}
On a les d\'eveloppements de Laurent suivants en $x=0$
$$
{\theta'_{\tau}\over\theta_{\tau}}(x) = {1\over x} - G_{2}(\tau)x 
- G_{4}(\tau)x^{3}- \cdots, 
\quad \wp_{\tau}(x) = {1\over x^{2}} + 3G_{4}(\tau)x^{2}+5G_{6}(\tau)x^{4}
+\cdots. $$
\end{lemma}

{\em D\'emonstration.} Le deuxi\`eme d\'eveloppement provient de 
$(x+a)^{-2} = a^{-2}-2xa^{-3}+\cdots$. 

D'apr\`es \cite{Po}, Thm. 3.9, $\wp_\tau = - (\sigma'/\sigma)'$, o\`u on pose $\sigma(z)
:= e^{{1\over 2}G_{2}(\tau)z^{2}}\theta_\tau(z)$. Donc $(\theta'_\tau/\theta_\tau)'(z) 
= -\wp_\tau(z)-G_{2}(\tau)$, 
ce qui d\'etermine le d\'eveloppement de $\theta'_\tau/\theta_\tau$ \`a une constante 
additive pr\`es. Cette constante est d\'etermin\'ee par le fait que  $\theta'_\tau/\theta_\tau$
est une fonction impaire. \hfill \qed\medskip  

On en d\'eduit 
$$
\tilde\wp_\tau(x)=\sum_{n\geq -1} (2n+1)G_{2n+2}(\tau)x^{2n}
=-({\theta'_\tau\over\theta_\tau})'(x), 
$$
o\`u on a pos\'e $G_{0}(\tau):= -1$. L'\'equation (\ref{formule:g:i}) implique alors
\begin{equation}\label{id:g:i}
g_i(0,\tau)=g_i(1,\tau)=-{1\over{2\pi\mathrm{i}}}\tilde\wp_\tau(x_i). 
\end{equation}

\begin{lemma}\label{newlemma}
On a l'identit\'e 
\begin{equation} \label{id:surv}
\forall x,y\in\CC, \quad (\partial_{x}\sigma^\tau_{x}) \sigma^\tau_{y} 
- \sigma^\tau_{x}(\partial_{y}\sigma^\tau_{y}) = \sigma^\tau_{x+y}(\wp_\tau(y)-\wp_\tau(x)) ,  
\end{equation}
\end{lemma}

{\em D\'emonstration.} Le membre de gauche a le m\^{e}me comportement que
le membre de droite sous les transformations de la variable muette $z\mapsto 
z+1$, $z\mapsto z+\tau$ ; pour \'etudier son comportement en $z=0$, on 
le transforme ainsi 
$$
(\partial_{x}\sigma_{x}^\tau) \sigma_{y}^\tau(z) 
- \sigma_{x}^\tau(\partial_{y}\sigma_{y}^\tau)(z) =
\sigma_{x}^\tau\sigma_{y}^\tau({{\partial_{x}\sigma_{x}^\tau}
\over{\sigma_{x}^\tau}} - {{\partial_{y}\sigma_{y}^\tau}
\over{\sigma_{y}^\tau}})(z)
= \sigma_{x}^\tau\sigma_{y}^\tau({\theta'_\tau\over\theta_\tau}(z+x)
- {\theta'_\tau\over\theta_\tau}(x) - {\theta'_\tau\over\theta_\tau}(z+y)
+ {\theta'_\tau\over\theta_\tau}(y))$$ 
dont le d\'eveloppement est 
${1\over {z^{2}}}\times z\times (({\theta'_\tau\over\theta_\tau})'(x)
- ({\theta'_\tau\over\theta_\tau})'(y)) + O(1) = 
{1\over z}(\wp_\tau(y)-\wp_\tau(x)) + O(1)$. On a donc un p\^{o}le simple en 
$0$, ce qui implique que le membre de gauche est proportionnel \`a 
$\sigma_{x+y}^\tau$ ; le d\'eveloppement en $0$ permet aussi de calculer le 
coefficient de proportionnalit\'e. \hfill\qed\medskip 

L'\'equation (\ref{id:g:i}) et le lemme \ref{newlemma} impliquent alors l'\'egalit\'e
$$
(g_i\omega_{i+1}-g_{i+1}\omega_i)(\tau,z)=(g_i(0,\tau)-g_{i+1}(0,\tau))\psi_{i,i+1}(\tau,z).
$$

On a alors: 

\begin{proposition} 
Les $(\omega_i)_{i=1,\ldots,n}$, $(g_i)_{i=1,\ldots,n}$ et $(\psi_{i,i+1})_{i=1,\ldots,n-1}$ de la d\'efinition \ref{def:ell} satisfont 
les hypoth\`eses de la proposition \ref{prop:3:1}. 
\end{proposition}

En appliquant la proposition \ref{prop:3:1}, on obtient la premi\`ere partie du r\'esultat suivant :

\begin{theorem} \label{thm:ode} 
\begin{equation} \label{I:ode}
(2\pi\on{i})\partial_{\tau} I_{x_{1},\ldots,x_{n}}(\tau)
= \tilde\wp_\tau(x_{1})I_{x_{2},\ldots,x_{n}}(\tau)
- \tilde\wp_\tau(x_{n})I_{x_{1},\ldots,x_{n-1}} (\tau)
+ \sum_{i=1}^{n-1}(\wp_\tau(x_{i+1})-\wp_\tau(x_{i}))I_{x_{1},\ldots,
x_{i}+x_{i+1},\ldots,x_{n}}(\tau). 
\end{equation}
et 
\begin{align*}
(2\pi\on{i})\partial_{\tau} J_{x_{1},\ldots,x_{n}}(\tau)
& = \tilde\wp_\tau(x_{1})J_{x_{2},\ldots,x_{n}}(\tau)
- \tilde\wp_\tau(x_{n})J_{x_{1},\ldots,x_{n-1}} (\tau)
+ \sum_{i=1}^{n-1}(\wp_\tau(x_{i+1})-\wp_\tau(x_{i}))J_{x_{1},\ldots,
x_{i}+x_{i+1},\ldots,x_{n}}(\tau)
\\
 & -{{2\pi\on{i}}\over\tau}(x_{1}\partial_{x_{1}}+\cdots + x_{n}
 \partial_{x_{n}})J_{x_{1},\ldots,x_{n}}(\tau). 
\end{align*}
\end{theorem}

La deuxi\`eme identit\'e est une cons\'equence de la premi\`ere identit\'e et de l'identit\'e 
modulaire (\ref{id:mod}), compte tenu de la relation modulaire 
$\tilde\wp_{-1/\tau}(x) = \tau^{2}\tilde\wp_{\tau}(\tau x) 
- 2\pi\on{i}\tau$, qui provient de $\wp_{-1/\tau}(x) = \tau^{2}
\wp_{\tau}(\tau x)$ et de $G_{2}(-1/\tau) = \tau^{2}G_{2}(\tau) - 2\pi\on{i}
\tau$ (\cite{Se}, \'equation (45) p. 156). \qed\medskip

\section{R\'ealisation de $\langle\delta_{2n},n\geq -1\rangle$ et comparaison de
syst\`emes diff\'erentiels}
\label{sect:real:fonct}

Comme $F = U(\f_{2}\ominus\CC x)\subset U(\f_{2})$ 
est un sous-$\f_{2}$-module de $U(\f_{2})$, on dispose d'un sous-espace 
$\on{Der}_{t}(\f_{2},F) \subset \on{Der}_{t}(\f_{2},U(\f_{2})) = 
\on{Der}_{t}(U(\f_{2}))$, qui est en fait une sous-alg\`ebre de Lie. 
($\on{Der}_{t}$ est l'ensemble des d\'erivations qui envoient $t= 
-[x,y]\in\f_{2}$ sur $0$.)  
D'autre part, le degr\'e en $y$ induit une graduation de ces alg\`ebres de Lie ; 
on note $\on{Der}(\f_{2})_{+}\subset \on{Der}(\f_{2})$ la partie
de $y$-degr\'e $>0$. 

On a donc une suite d'inclusions d'alg\`ebres de Lie 
$$
\langle\delta_{2n},n\geq -1\rangle \subset \on{Der}_{t}(\f_{2})_{+}
\subset \on{Der}_{t}(\f_{2},F) \subset \on{Der}_{t}(F). 
$$
Le but de cette section est d'\'etablir un isomorphisme
$$
\on{Der}_{t}(\f_{2},F)\simeq {\mathcal G}_{0}
$$
entre $\on{Der}_{t}(\f_{2},F)$ et une alg\`ebre de Lie ``fonctionnelle'' 
${\mathcal G}_{0}$ explicite, puis d'en d\'eduire le lien entre les 
\'equations diff\'erentielles du th\'eor\`eme \ref{thm:ode} et celles satisfaites
par $A(\tau)$, $B(\tau)$ (\'equations (\ref{ED:A:B})). 

Pour $n\geq 1$, on pose ${\mathcal G}[n]:=\CC(x_{1},\ldots,x_{n+1})$, le corps des fractions rationnelles 
\`a coefficients dans $\CC$ en $n+1$ ind\'etermin\'ees $x_1,\ldots,x_{n+1}$. On pose 
$$
{\mathcal G}:= \oplus_{n\geq 1}{\mathcal G}[n]. 
$$

L'\'enonc\'e suivant est imm\'ediat. 
\begin{proposition}
${\mathcal G}$ est munie d'une 
structure d'alg\`ebre de Lie gradu\'ee donn\'ee par
$$
[\varphi,\psi]:= 
\sum_{i=1}^{m+1}\varphi^{i,i+1,\ldots,i+n}
\psi^{1,\ldots,i-1,ii+1\ldots i+n,i+n+1,\ldots,n+m+1}
- ((\varphi,n)\leftrightarrow(\psi,m))\in{\mathcal G}[n+m]    
$$
pour $\varphi\in{\mathcal G}[n]$, $\psi\in{\mathcal G}[m]$ ;  
on note $\varphi^{1,\ldots,n+1}:= \varphi(x_{1},\ldots,x_{n+1})$, 
$\varphi^{12,3}:= \varphi(x_{1}+x_{2},x_{3})$, etc. 
\end{proposition}

L'espace $F_{\infty}:= \CC(x_{i},i\in\ZZ)$ des fractions rationnelles en une 
infinit\'e de variables est un ${\mathcal G}$-module via 
$\varphi * f:= \sum_{i\in\ZZ} \varphi^{i,i+1,\ldots,i+n}
f^{\ldots,i-1,ii+1\ldots i+n,i+n+1,\ldots}$.  

Par l'identification de $\CC[x_{1},\ldots,x_{n}]$ \`a $\CC[x_{1},\ldots,x_{n+1}]/
(x_{1}+\cdots+x_{n+1})$, on obtient une action du groupe sym\'etrique $S_{n+1}$ 
sur cette premi\`ere alg\`ebre (en d'autres termes, on a affaire \`a l'alg\`ebre sym\'etrique
du quotient $\CC^{n+1}/\CC$ de la repr\'esentation naturelle par la triviale). 
On note $C_{n+1}\subset S_{n+1}$ le sous-groupe cyclique. L'\'el\'ement
$x_1\cdots x_n(x_1+\cdots+x_n)$ de $\CC[x_1,\ldots,x_n]$ est invariant par l'action de ce groupe. 

Pour $n\geq 1$, on pose 
$$
{\mathcal G}_0[n]:= \big((x_{1}\cdots x_{n}(x_{1}+\cdots+x_{n}))^{-1}\CC[x_{1},\ldots,x_{n}]\big)^{C_{n+1}},  
$$
o\`u l'espace entre parenth\`eses est celui des fractions rationnelles en $x_1,\ldots,x_n$
avec d\'enominateur $x_1\cdots x_n(x_1+\cdots+x_n)$. 
On pose alors 
$$
{\mathcal G}_{0}:= \oplus_{n\geq 1}{\mathcal G}_{0}[n]. 
$$

\begin{proposition}
Une structure d'alg\`ebre de Lie gradu\'ee est d\'efinie sur ${\mathcal G}_{0}$
par
\begin{align*}
[\varphi,\psi]_{0} & := 
\sum_{i=1}^{n}(\varphi^{i,i+1,\ldots,i+n-1} - \varphi^{i+1,\ldots,i+n})
\psi^{1,\ldots,i-1,ii+1\ldots i+n,i+n+1,\ldots,n+m} 
\\
 & - \sum_{j=1}^{m}(\psi^{j,j+1,\ldots,j+m-1} - \psi^{j+1,\ldots,j+m})
\varphi^{1,\ldots,j-1,jj+1\ldots j+m,j+m+1,\ldots,n+m} 
 \\
 & - \varphi^{1,\ldots,n}\psi^{n+1,\ldots,n+m}
+\varphi^{m+1,\ldots,n+m}\psi^{1,\ldots,m}\in{\mathcal G}_0[n+m]
\end{align*}
pour $\varphi\in{\mathcal G}_0[n]$, $\psi\in{\mathcal G}_0[m]$.  
L'espace $F$ a une structure de ${\mathcal G}_{0}$-module gradu\'e par 
\begin{align} \label{act:F}
\varphi \bullet f & \nonumber := 
\sum_{i=1}^{m}(\varphi^{i,i+1,\ldots,i+n-1} 
- \varphi^{i+1,\ldots,i+n})f^{1,\ldots,i-1,ii+1\ldots i+n,i+n+1,\ldots,n+m} \\
 & - \varphi^{1,\ldots,n}f^{n+1,\ldots,n+m} + \varphi^{n+1,\ldots,n+m}
f^{1,\ldots,m}\in F_{n+m} 
\end{align}
pour $\varphi\in{\mathcal G}_0[n]$, $f\in F_n$. 
Un isomorphisme d'alg\`ebres de Lie gradu\'ees $\on{Der}_{t}(\f_{2},F)\simeq 
{\mathcal G}_{0}$ est donn\'e par 
$$
\on{Der}_{t}(\f_{2},F)[n]\ni D
\leftrightarrow \varphi(x_{1},\ldots,x_{n})
\in{\mathcal G}_{0}[n], 
$$
o\`u $D$ est la d\'erivation donn\'ee par $x\mapsto u$, $y\mapsto v$, 
o\`u $u\in F_n$, $v\in F_{n+1}$ sont donn\'es par 
$$
u(x_{1},\ldots,x_{n}) = (x_{1}+\cdots+x_{n})\varphi(x_{1},\ldots,x_{n}), 
$$
$$
v(x_{1},\ldots,x_{n+1}) = ({1\over {x_{1}}} - {1\over{x_{1}+\cdots
+x_{n+1}}})\varphi(x_{2},\ldots,x_{n+1})
+ ({1\over{x_{1}+\cdots
+x_{n+1}}} - {1\over {x_{n+1}}})\varphi(x_{1},\ldots,x_{n}). 
$$
\end{proposition}

Un morphisme d'alg\`ebres de Lie ${\mathcal G}_{0}\to{\mathcal G}$
est par ailleurs donn\'e par
$$
{\mathcal G}_{0}[n]\ni \varphi(x_{1},\ldots,x_{n})\mapsto 
\varphi(x_{1},\ldots,x_{n}) - \varphi(x_{2},\ldots,x_{n+1})
\in{\mathcal G}[n]. 
$$

{\em D\'emonstration.} Soit $D \in\on{Der}_{t}(\f_{2},F)[n]$. Cet \'el\'ement est d\'etermin\'e 
par le couple $(u,v):= (D(x),D(y))\in F_{n}\times F_{n+1}$. La condition sur $(u,v)$ est 
$$
(x_{1}+\cdots+x_{n+1})v(x_{1},\ldots,x_{n+1}) = x_{1}^{-1}u(x_{2},\ldots,
x_{n+1}) - x_{n+1}^{-1}u(x_{1},\ldots,x_{n}) 
$$
(identit\'e dans $(x_{1}\cdots x_{n+1})^{-1}\CC[x_{1},\ldots,x_{n+1}]$). 
On dispose d'une application de r\'eduction mo\-dulo $x_{1}+\cdots+x_{n+1}$
de cet espace vers $(x_{1}\cdots x_{n}(x_{1}+\cdots+x_{n}))^{-1}\CC[x_{1},
\ldots,x_{n}]$. L'image de cette identit\'e exprime alors la $C_{n+1}$-invariance
de $\varphi(x_{1},\ldots,x_{n}):= 
u(x_{1},\ldots,x_{n})/(x_{1}+\cdots+x_{n})$. On a donc une application 
lin\'eaire 
\begin{equation} \label{appli}
\on{Der}_{t}(\f_{2},F)[n] \to 
(x_{1}\cdots x_{n}(x_{1}+\cdots+x_{n}))^{-1}\CC[x_{1},\ldots,x_{n}]^{C_{n+1}},
\end{equation} 
$D\mapsto \varphi$. Cette application est 
injective car la nullit\'e de $u$ implique celle de $v$. Les deux derni\`eres formules  
de la proposition d\'efinissent une application 
$$
(x_{1}\cdots x_{n}(x_{1}+\cdots+x_{n}))^{-1}\CC[x_{1},\ldots,x_{n}]^{C_{n+1}}
\to F_{n}\times F_{n+1}
$$ 
(le p\^{o}le en $x_{1}+\cdots+x_{n+1}$ disparaissant 
par $C_{n+1}$-invariance), qui est en fait d'image dans $\on{Der}_{t}(\f_{2},F)[n]$
et inverse \`a (\ref{appli}). On v\'erifie alors que le transport \`a 
${\mathcal G}_{0}$ de la structure d'alg\`ebre de Lie sur $\on{Der}_{t}(\f_{2},F)$ 
et de module de $F$ sur cette alg\`ebre de Lie est donn\'e par les formules de l'\'enonc\'e. 
\hfill \qed\medskip 

On a 
$$
{\mathcal G}_{0}[1] = x_{1}^{-2}\CC[x_{1}^{2}]. 
$$

\begin{lemma}
L'isomorphisme ${\mathcal G}_{0}\simeq \on{Der}_{t}(\f_{2},F)$ induit 
la correspondance 
$$
\delta_{2n}\in \on{Der}_{t}(\f_{2})_{+}[1] \subset
\on{Der}_{t}(\f_{2},F)[1] \leftrightarrow {\mathcal G}_{0}[1]
\ni x_{1}^{2n}. 
$$
\end{lemma}

{\em D\'emonstration.} La d\'erivation correspondant \`a $x_{1}^{2n}$ est une
d\'erivation de $\on{Der}_{t}(\f_{2},F)$ telle que $x\mapsto u = 
x_{1}^{2n+1} \leftrightarrow [x^{2n+2}y]$. Comme l'application 
$\on{Der}_{t}(\f_{2},F)\to F$, $D\mapsto D(x)$
est injective, cette d\'erivation coincide avec la d\'erivation $\delta_{2n}$ d\'efinie en section \ref{sect:ED}. 
\hfill \qed\medskip 

Rappelons par ailleurs la correspondance 
$$
e^{\on{i}\pi t}A(\tau)\in U(\f_{2}\ominus \CC x)\leftrightarrow 
F\ni ((-1)^{n}I_{x_{n},\ldots,x_{1}}(\tau))_{n\geq 0} =: \tilde I(\tau). 
$$
(section \ref{sect:corr}). 
D'apr\`es (\ref{ED:A:B}) et l'invariance de $t$ sous les $\delta_{2n}$, 
$n\geq -1$, $e^{\on{i}\pi t}A(\tau)$ satisfait l'\'equation diff\'erentielle
$$
2\pi\on{i}\partial_{\tau}(e^{\on{i}\pi t}A(\tau)) = -(\sum_{n\geq -1} (2n+1) G_{2n+2}(\tau)
\delta_{2n}) (e^{\on{i}\pi t}A(\tau)). 
$$ 
L'image de cette \'equation diff\'erentielle sous l'isomorphisme 
$U(\f_{2}\ominus \CC x)\simeq F$ donne
$$
2\pi\on{i}\partial_{\tau}\tilde I(\tau) = -(\sum_{n\geq -1} (2n+1) G_{2n+2}(\tau)
x_{1}^{2n}) \bullet \tilde I(\tau) = -\tilde\wp_{\tau}(x_{1}) \bullet \tilde I(\tau),  
$$
donc si $I(\tau):= (I_{x_{1},\ldots,x_{n}}(\tau))_{n\geq 0}$, alors
$2\pi\on{i}\partial_{\tau}I(\tau) = -\tilde\wp_{\tau}(x_{1})\bullet
I(\tau)$, c'est-\`a-dire que pour chaque $n$
$$
2\pi\on{i}\partial_{\tau}I_{x_{1},\ldots,x_{n}}(\tau) 
=  -\tilde\wp_{\tau}(x_{1}) 
\bullet I_{x_{1},\ldots,x_{n-1}}(\tau). 
$$
Compte tenu de la formule (\ref{act:F}) pour l'action de ${\mathcal G}_{0}$ 
sur $F$, on retrouve ainsi 
la premi\`ere \'equation diff\'erentielle du th\'eor\`eme \ref{thm:ode}. 

De m\^{e}me, $e^{-\on{i}\pi t}B(\tau)$ satisfait l'\'equation diff\'erentielle
$$
2\pi\on{i}\partial_{\tau}(e^{-\on{i}\pi t}B(\tau)) = -(\sum_{n\geq 1}
(2n+1)G_{2n+2}(\tau)\delta_{2n})(e^{-\on{i}\pi t}B(\tau)). 
$$
Soit $\tilde B(\tau):= \on{exp}({2\pi\on{i}\over\tau}e_{+})(e^{-\on{i}\pi t}
B(\tau))$ (voir section \ref{sect:corr}), on en d\'eduit 
$$
2\pi\on{i}\partial_{\tau}\tilde B(\tau) = 
-({2\pi\on{i}\over\tau}h + \sum_{n\geq -1}(2n+1)G_{2n+2}(\tau)\delta_{2n})
(\tilde B(\tau)), $$
o\`u $h := [e_{+},\delta]$ est la d\'erivation de $\f_{2}$ donn\'ee par 
$(x,y)\mapsto (x,-y)$, compte tenu de $[e_{+},\delta_{2n}]=0$ si $n\geq 0$
et ${1\over 2}[e_{+},[e_{+},\delta_{-2}]] + e_{+}=0$. 

On a la correspondance 
$$
U(\f_{2} \ominus \CC x)\ni \tilde B(\tau)\leftrightarrow ((-1)^{n}J_{x_{n},\ldots,x_{1}}(\tau))=:\tilde 
J(\tau)\in F, 
$$
par ailleurs la d\'erivation $h$ se transporte sous cette correspondance en 
la d\'erivation de $F = \oplus_{n\geq 0}F_{n}$ de degr\'e z\'ero, 
op\'erant sur $F_{n}$ comme $\xi:= \sum_{i=1}^{n} x_{i}\partial_{x_{i}}$. 
On en d\'eduit 
$$
2\pi\on{i}\partial_{\tau}\tilde J(\tau) = -({2\pi\on{i}\over\tau}\xi
+\tilde\wp_{\tau}(x_{1})\bullet)\tilde J(\tau), 
$$
donc $J(\tau):= (J_{x_{1},\ldots,x_{n}}(\tau))_{n\geq 0}$ satisfait la m\^{e}me
\'equation diff\'erentielle, donc 
$$
2\pi\on{i}\partial_{\tau} J_{x_{1},\ldots,x_{n}}(\tau) = 
- {2\pi\on{i}\over\tau}(\sum_{i=1}^{n}x_{i}\partial_{x_{i}})
J_{x_{1},\ldots,x_{n}}(\tau)
-\tilde\wp_{\tau}(x_{1})\bullet J_{x_{1},\ldots,x_{n-1}}(\tau),  
$$
ce qui permet de retrouver la deuxi\`eme \'equation diff\'erentielle du th\'eor\`eme
\ref{thm:ode}.  

\section{D\'eveloppement asymptotique des analogues elliptiques des nombres 
multiz\'etas}
\label{section:DA}

Dans cette section, nous utilisons les \'equations diff\'erentielles
satisfaites par les fonctions $A(\tau)$ et $B(\tau)$ (\'equations (\ref{ED:A:B})) et 
leur comportement \`a l'infini ((\ref{comp:A}), (\ref{comp:B}))
pour en obtenir un d\'eveloppement asymptotique en $\tau\to\on{i}\infty$. 
Nous en d\'eduisons la forme du d\'eveloppement asymptotique des 
fonctions $I_{\ul{d}}(\tau)$, $J_{\ul{d}}(\tau)$ dans cette r\'egion. 

\subsection{D\'eveloppement de $g(\tau)$}  \label{DA:g}

Soit ${\mathfrak G}$ la compl\'etion de l'alg\`ebre de Lie 
$\langle\delta_{2n},n\geq -1\rangle
\subset \on{Der}_{t}(\f_{2})$ pour le bidegr\'e en $(x,y)$ ; on a 
$|\delta_{2n}| = (2n+1,1)$. Soit $G:= \on{exp}({\mathfrak G})
\subset \on{Aut}_{t}(\hat\f_{2})$ le groupe de Lie correspondant. 

\begin{proposition} \label{prop:h}
Il existe une unique fonction $g(\tau) : \HH\to G$, telle que 
$$
2\pi\on{i}\partial_{\tau}g(\tau) = -(\sum_{n\geq -1}
(2n+1)G_{2n+2}(\tau)\delta_{2n})g(\tau)$$
et $g(\tau)\simeq e^{{-1\over{2\pi\on{i}}}(\delta_{-2}+\sum_{n\geq 0}
(2n+1)\cdot 2\zeta(2n+2)\delta_{2n})\tau}= e^{D_{0}\tau}$ en $\tau\to
\on{i}\infty$. Il existe une collection $(h_{k})_{k\geq 0}$, avec $h_{0}=1$, 
telle que $g(\tau)$ a le d\'eveloppement asymptotique 
$$
g(\tau)\simeq \sum_{k,n\geq 0}{1\over{n!}}h_{k}D_{0}^{n}\tau^{n}e^{2\pi\on{i}
k\tau}
$$
en $\tau\to\on{i}\infty$. 
\end{proposition}

{\em D\'emonstration.} Posons $D(\tau):= {{-1}\over{2\pi\on{i}}}
\sum_{n\geq -1}(2n+1)G_{2n+2}(\tau)\delta_{2n}$. Posons si $m\geq 1$, 
$g_{2m}(n):= {{2(2\pi\on{i})^{2m}}\over{(2m-1)!}}\sigma_{2m-1}(n)$
(o\`u $\sigma_{k}(n) = \sum_{d|n}d^{k}$) si $n>0$, et $g_{2m}(0):= 
2\zeta(2m)$ ; et posons $g_{0}(n)=0$ si $n>0$, et $g_{0}(0)=-1$. 
Alors $G_{2m}(\tau) = \sum_{n\geq 0} g_{2m}(n)e^{2\pi\on{i}n\tau}$ et 
$$D(\tau) = \sum_{m\geq 0}D_{m}e^{2\pi\on{i}m\tau}, 
\text{\ o\`u\ }D_{m}:= 
{{-1}\over{2\pi\on{i}}} \sum_{n\geq -1}
(2n+1)g_{2n+2}(m)\delta_{2n}.
$$
Comme $D_{0}$ est de $y$-degr\'e $1$, $2\pi\on{i}m-\on{ad}D_{0}$
est inversible dans $\on{End}(U{\mathfrak G})$ si $m>0$. D\'efinissons 
$(h_{m})_{m\geq 0}$ par $h_{0}:=1$, 
$$
h_{m}:= (2\pi\on{i}m - \on{ad}D_{0})^{-1}(\sum_{m'+m''=m\atop{m'>0}}
D_{m'}h_{m''})\text{\ si\ }m>0.
$$
Alors $h(\tau):= \sum_{m\geq 0}h_{m}e^{2\pi\on{i}m\tau}$ est une solution 
formelle de $\partial_{\tau}h(\tau) = D(\tau)h(\tau) - h(\tau)D_{0}$, qui est 
l'\'equation diff\'erentielle satisfaite par $g(\tau)e^{-D_{0}\tau}$ ; cette 
fonction admet donc $h(\tau)$ comme d\'eveloppement asymptotique. \hfill 
\qed\medskip 

\subsection{D\'eveloppements de $A(\tau)$, $B(\tau)$} \label{sect:DA:A:B}

Posons 
$$
A_{\infty}:= \Phi(\tilde y,t)e^{2\pi\on{i}\tilde y}\Phi(\tilde y,t)^{-1}, 
\quad
\underline B(\tau):= e^{\on{i}\pi t}\Phi(-\tilde y-t,t)e^{2\pi\on{i}x}
e^{2\pi\on{i}\tilde y\tau}\Phi(\tilde y,t)^{-1}, 
\quad B_{\infty}:= B[0]. 
$$
D'apr\`es \cite{CEE},  on a 
$$
A_{\infty} = e^{\tau D_{0}}(A_{\infty}), \quad 
\underline B(\tau) = e^{\tau D_{0}}(B_{\infty}). 
$$
Posons $h(\tau):= g(\tau)e^{-\tau D_{0}}$, alors d'apr\`es
la proposition \ref{prop:h}, $h(\tau)$ admet le d\'eveloppement asymptotique 
$h(\tau)\simeq 1+\sum_{m>0}h_{m}e^{2\pi\on{i}m\tau}$. 

On a alors 
$A(\tau) = g(\tau)(A_{\infty}) = h(\tau)(A_{\infty})$
donc $A(\tau)$ admet le d\'eveloppement asymptotique
$$ \label{dev:A}
A(\tau)\simeq \sum_{m\geq 0}e^{2\pi\on{i}m\tau} h_{m}(A_{\infty}), 
$$
et $B(\tau) = g(\tau)(B_{\infty}) = h(\tau)(\underline B(\tau))$, donc
$\tilde B(\tau)$ admet le d\'eveloppement asymptotique  
$$ \label{dev:B}
\tilde B(\tau) \simeq \on{exp}(-{{2\pi\on{i}}\over\tau}e_{+})h(\tau)
(\underline B(\tau)). 
$$
 
\subsection{D\'eveloppements de $I_{\ul{d}}(\tau)$, $J_{\ul{d}}(\tau)$} 
\label{sect:DA:I:J}

Soit ${\bold k}_{MZV}\subset \CC$ le $\QQ$-sous-anneau engendr\'e par les 
multiz\'etas. L'associateur $\Phi$ \'etant \`a coefficients dans ${\bold k}_{MZV}$, on 
d\'eduit de (\ref{dev:A}) et (\ref{dev:B}) : 

\begin{proposition}
Les fonctions $I_{\ul{d}}(\tau)$, $J_{\ul{d}}(\tau)$ admettent les d\'eveloppements
asymptotiques 
$$
I_{\ul{d}}(\tau)\simeq \sum_{n\geq 0} I_{\ul{d},n}e^{2\pi\on{i}n\tau}, 
\quad 
J_{\ul{d}}(\tau) \simeq \sum_{n\geq 0}\sum_{s\in\ZZ}J_{\ul d,n,s}
\tau^{s}e^{2\pi\on{i}n\tau} ,  $$
dans lesquels les coefficients sont dans ${\bold k}_{MZV}[2\pi\on{i}]$. 
Dans la deuxi\`eme s\'erie, la deuxi\`eme somme $\sum_{s}$ est finie pour 
tout $n\geq 0$. 
\end{proposition}

\def\refname{Bibliographie}
\end{document}